\documentclass[12pt]{amsart}
\usepackage{amsthm}
\usepackage{epsfig}
\usepackage{amssymb}
\usepackage{latexsym}
\usepackage{amsmath}
\newtheorem{theorem}{Theorem} 
\newtheorem{definition}[theorem]{Definition}

\newtheorem{corollary}[theorem]{Corollary}
\newtheorem{lemma}[theorem]{Lemma}

\newtheorem{proposition}[theorem]{Proposition}

\def \vp{\varphi}

\def\R{{\Bbb R}}

\def\C{{\Bbb C}}

\def\Z{{\Bbb Z}}
\def\P{{\Bbb P}}

\def \ta{\tau}
\def \ta1{\tau_1}

\def \dl{\delta}
\def \g{\gamma}
\def \G{\Gamma}

\def \vp{\varphi}

\def\G{\Gamma}

\def \la{\langle}
\def \ra{\rangle}

\def \ra{\rangle}

\def \tC{\tilde{C}}

\def \tP{\tilde{\pi}}

\def \A{{\mathcal A}}

\def \prodl{\prod\limits}
\def\hF{\hat{F}}
\def\uZ{\underline{Z}}
\def\bZ{\bar{Z}}

\def\uZut{\underline{Z}^2}
\def\uZumt{\underline{Z}^{-2}}

\def \zovera {
    \mathop{\lower 10pt \hbox{${\buildrel{\displaystyle\bar{z}} \over {\scriptstyle{(a)}}} $}}
    {\lower 4pt \hbox{${\scriptstyle{ij}}$}} 
} 

\def\vK{{van Kampen}}
\def\sub{{\subseteq}}

\newcommand\set[1]{{\{{#1}\}}}
\def\suchthat{{\,:\,}}

\def\wrt{{with respect to }}

\def\st{{such that }}
\def\pitil{{\tilde{\pi}_1}}

\newcommand\sg[1]{{\left<{#1}\right>}}

\newcommand\begintable[1][] {{}}

\newif\ifXY 
\XYtrue     
\newif\ifbigmatrices
\bigmatricestrue 

 \ifXY
 \usepackage{xy}
 \fi
 \ifXY
 \xyoption{all}
 \fi

\begin{document}

\title {The fundamental group of the complement of the branch curve of $T \times T$ in $\C^2$}
 \author{ Meirav Amram \and Mina Teicher}
\address{Meirav Amram, Mathematisches Institut, Bismarck Strasse 1  1/2, Erlangen, Germany}
\email{meirav@macs.biu.ac.il, amram@mi.uni-erlangen.de}
\address{Mina Teicher, Mathematics department, Bar-Ilan university, 52900 Ramat-Gan, Israel}
\email{teicher@macs.biu.ac.il}

\renewcommand{\subjclassname}{%
       \textup{2000} Mathematics Subject Classification}


\date{\today}

\maketitle

\begin{abstract}
This paper is the second in a series of three papers concerning the surface $T \times T$, where $T$ is a complex torus.
We compute the fundamental group of the branch curve of the surface in $\C^2$, using the van Kampen Theorem and the braid monodromy factorization of the curve.
\end{abstract}

\section{Background}
The concentration on the fundamental group of a complement of a
branch curve of an algebraic surface $X$ with respect to a generic projection 
onto $\C\P^2$,
leads us to the computation of $\pi_1(X_{Gal})$ the fundamental
group of the Galois cover of 
$X$ with respect to this generic projection. 
Galois covers are surfaces of a general type.

Bogomolov conjectured that the Galois covers corresponding 
to generic projections of algebraic surfaces to $\C \P^2$ have infinite fundamental groups. 

In [7] we justify Bogomolov's conjecture by proving 
that  $\pi_1(T \times T)_{Gal}$ is an infinite group. 
 
In order to compute $\pi_1(T \times T)_{Gal}$, we have to enclode the braid monodromy factorization of the branch curve $S$ of $T \times T$. Then we have to apply the van Kampen Theorem on the factors in the factorization in order to get  relations for $\pi_1(\C^2-S,*)$ 
the fundamental group of the complement of $S$ in $\C^2$.

The fundamental group of the Galois cover $X_{Gal}$ is known to 
be a quotient  of a certain subgroup of the fundamental
group of the complement of $S$.

We recall shortly the computations from [6].

Let $X=T \times T$ be an algebraic surface (where $T$ is a 
complex torus) embedded in $\C\P^5$, 
and $f:X\rightarrow
\C\P^2$ be a generic projection. 
We degenerate $X$ to a union of $18$ planes $X_0$([6, Section 3]). 
We numerate the lines and vertices as shown in Figure 1.
\begin{figure}[h]\label{nums}
\epsfxsize=5cm 
\epsfysize=5cm 
\begin{center}
\epsfbox {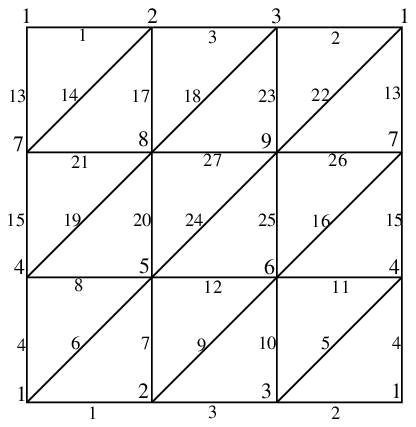}
\end{center}
\caption{}\end{figure}

We have a generic projection  $f_0:X_0\rightarrow\C\P^2$. 
We get a degenerated branch curve $S_0$ which is a line arrangement and 
compounds nine 6-points. We regenerate each 6-point separately.

We concentrate for example 
in a regeneration in a neighbourhood of $V_2$. 
We consider the local numeration of lines meeting at $V_2$ 
([6, Figure 6]). First,
the diagonal lines 4 and 5 become conics which are tangent 
to the lines 2,3 and 1,6 respectively, see Figure 2.

\begin{figure}[h]\label{1step}
\epsfxsize=11cm 
\epsfysize=9cm 
\begin{center}
\epsfbox {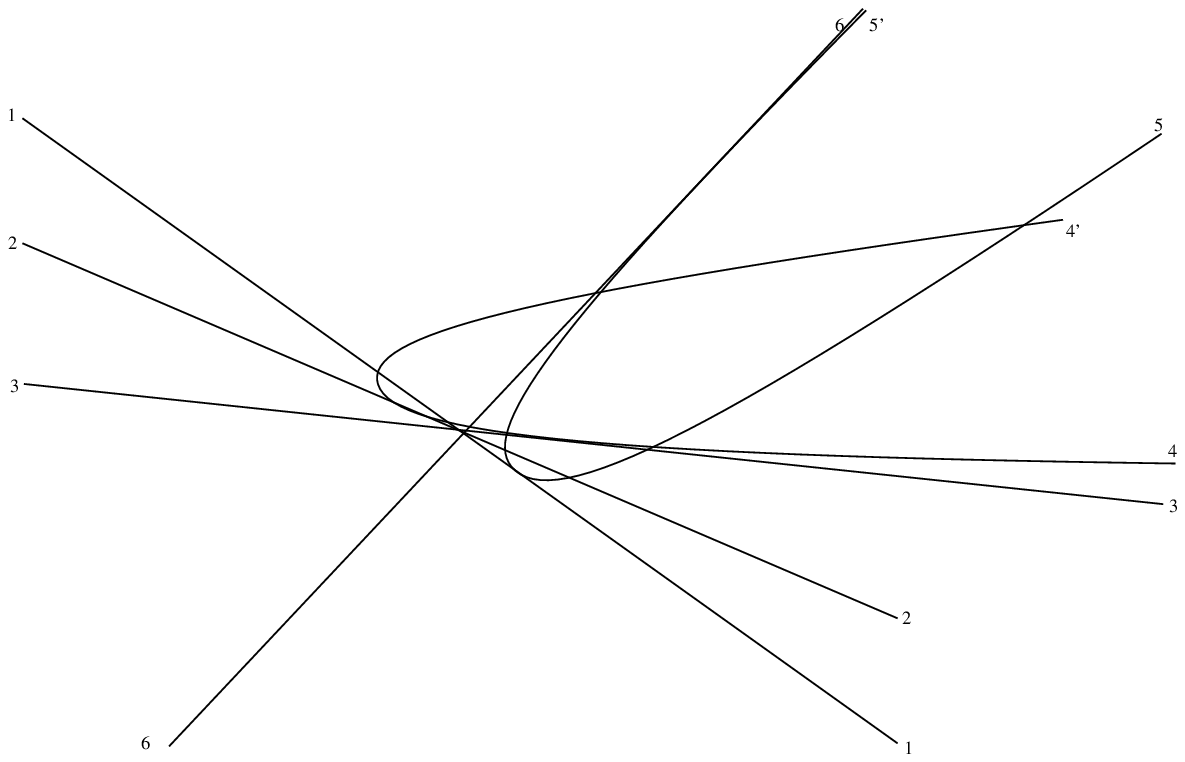}
\end{center}
\caption{}\end{figure}

Now, concentrating in a neighbourhood of the left 4-point (the
intersection point of 1,2,3,6), the two horizontal lines become a 
hyperbola.
Each one of the two vertical lines is replaced by two parallel lines, 
which are  tangent to the hyperbola, see Figure 3.

\begin{figure}[h]\label{2step}
\epsfxsize=11cm 
\epsfysize=9cm 
\begin{minipage}{\textwidth}
\begin{center}
\epsfbox {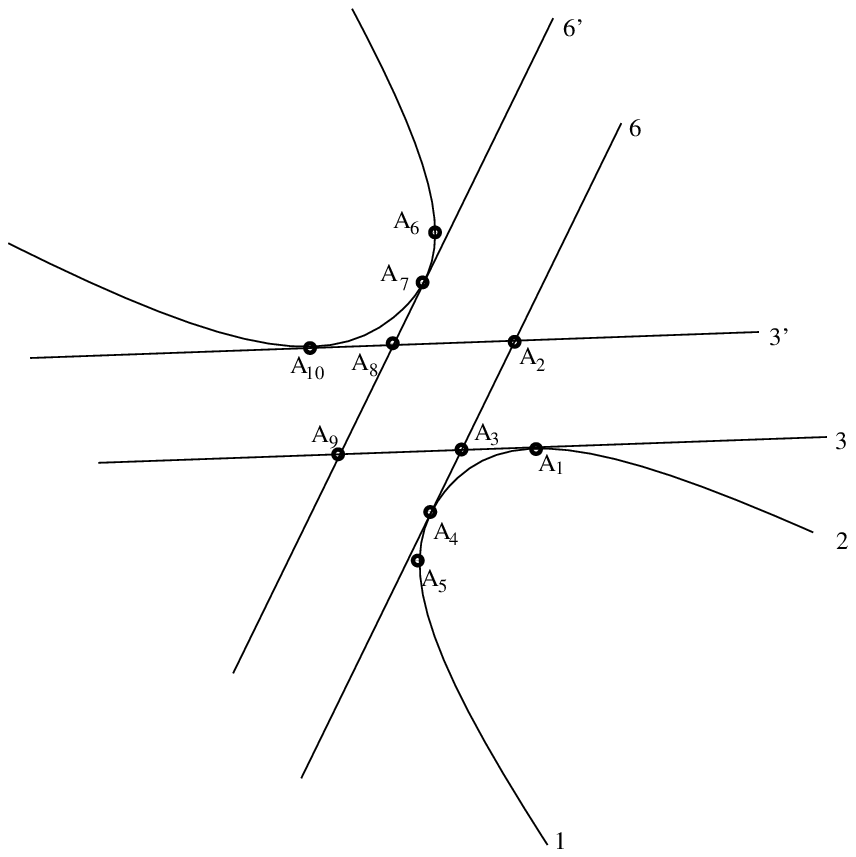}
\end{center}
\end{minipage}
\caption{}\end{figure}

Finally, the hyperbola is doubled.  
Each one of the tangent points is replaced by three cusps
([6, Theorem 13]).  Naturally, 
it occurs that each node is replaced by two or four ones.

We end up with a regenerated cuspidal curve , 
which has a degree of 12.

We do so to each one of the 6-points and get regenerated curves 
$S_i, 1 \leq i \leq 9$.
The union of these curves is the regenerated branch curve $S$.
Deg$S_0 = 27$, thus deg$S=54$.  The reason is that each intersection
point $q_j$ of the curve $S_0$ with the typical 
fiber was replaced by two close points $\{q_{j}, q_{j'}\}$.

In order to define a g-base for the fundamental group of the complement of $S$ in $\C^2$, we present the following situation(following Figure 4 to understand the below notions).\\
 $S$ is an algebraic curve in $\C^2\ , \ {54} = \mbox{deg } S$.
 $\pi: \C^2 \rightarrow \C$ a generic projection on the first coordinate.
 $K(x) = \{y \mid (x,y) \in S\}$ is the projection to the $y$-axis of $\pi^{-1}(x) \cap S$. Let
 $N = \{x \mid \# K(x) < {54} \}$ and
 $M' = \{ x \in S \mid \pi_{\mid x} \mbox{ is not \'{e}tale at } x \}$ such that $\pi (M') = N$.
 Assume $\# (\pi^{-1}(x) \cap M') =1, \forall x \in N$. Let $E$
 (resp. $D$) be a closed disk  on $x$-axis (resp. $y$-axis), \st
 $M' \subset E \times D, N \subset \mbox{ Int } (E)$. We choose $u
 \in \partial E, \ x \ll u \;\;\; \forall x \in N$. $\C_u = \{q_1,
q_{1'}, \cdots, q_{27}, q_{27'}\}$.

\begin{figure}[h]\label{S}
\epsfxsize=7cm 
\epsfysize=6cm 
\begin{center}
\epsfbox {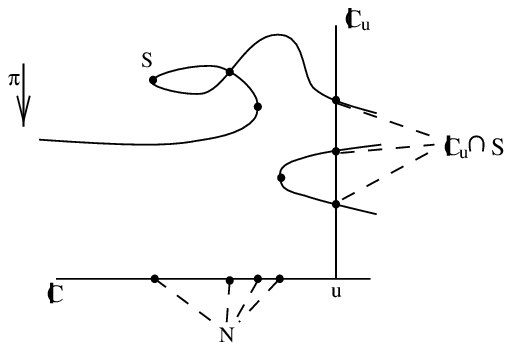}
\end{center}
\caption{}\end{figure}

We now specify a standard set of generators for
the fundamental group $\pi_1(\C^2-S,M)$, where $M$ is a point 
outside $S$.

Write $S \cap \C_u =\set{q_1, \cdots , q_{27'}}$. 
 Let $\g_j$ be paths from $M$ to $q_j$ $\forall j$, 
\st the $\g_j$ do not meet each other in any point except $M$.  
Let $\eta_j$ be a  small oriented circle around $q_j$.  
Let $\g'_j$ be the part of $\g_j$ outside $\eta_j$, and take
$\G_j = \g'_j \eta_j (\g'_j)^{-1}$. In the same way we specify 
generators  $\G_{j'}$.  
The set $\{\G_j,\G_{j'}\}_{j=1}^{27}$ freely generates 
$\pi_1(\C_u - S, M)$ [40].  Such a set is called a g-base 
for $\pi_1(\C_u - S, M)$.

By Lemma \ref{lm43} we have a surjection
$\pi_1(\C_u - S, M)
\stackrel{\nu}{\rightarrow} \pi_1(\C^2 - S, M) \rightarrow 0$, 
so the  set $\{\nu(\G_j)\}$ generates  $\pi_1(\C^2 - S,M)$. 
By abuse of notation, we shall denote $\nu(\G_j)$ by $\G_j$.
A presentation for $\pi_1(\C^2-S,M)$  is obtained by the
\vK\ Theorem, from a list of braids in $B_{54}[\C_u, \C_u \cap S]$.

\smallskip

The group $\pi_1(\C^2 - S,M)$ acts on the points in $\C_u$. This leads to  
a permutation representaion $\psi: \pi_1(\C^2-S,M) \rightarrow S_{18}$, 
$18$ is the number of planes in $X_0$.

Let $\sg{\G^2_j,\G^2_{j'} }$ denote the normal subgroup generated by $\G^2_j,\G^2_{j'} $.

\begin{definition}\label{pitil}

Define $$\pitil = {\pi_1(\C^2-S,M) \over {\sg{\Gamma_j^2, \G_{j'}}}}.$$

\end{definition}

 Since $f$ is stable, $\G_j$ and $\G^2_{j'}$ induce 
a transposition in $S_{18}$
so that $\sg{\G^2_j} \sub \ker \psi$.
The map $\pitil \rightarrow S_{18}$ is also denoted $\psi$.

 By the  isomorphism theorems, we have an exact sequence
 \begin{equation} \label{exact1}
 1
 \rightarrow \ker \psi
 \rightarrow \pitil
 \stackrel{\psi}{\rightarrow} S_{18} \rightarrow 1.
 \end{equation}

\begin{definition}\label{galois}
Consider the fibered product\\
$$\underbrace{X \times_f \cdots \times_f X}_{n}  = \set{(x_1, \cdots , x_n) 
\in X^n \suchthat f(x_1) = \cdots = f(x_n)},$$ and  the diagonal 
$$\Delta = \set{(x_1, \cdots , x_n) \in X \times_f \cdots \times_f X 
\mid x_i = x_j \mbox{ for some } i \neq j}.$$
The surface  $X_{Gal}$ is the Galois cover    
of $X$ with respect to the generic projection $f :X  \rightarrow \C\P^2$.   
That is the Zariski closure of the complement of $\Delta$: 
$$X_{Gal}=\overline{X \times_f \cdots\times_f X - \Delta}.$$
\end{definition}

Let $X_{Gal}^{Aff}$ be the part of $X_{Gal}$
 lying over $\C^2$ ($\subseteq \C \P^2$).  
 There is a surjective map $X_{Gal}^{Aff} \rightarrow X_{Gal}$.

\begin{theorem}{\bf{[16, Secion 0.3]}}\label{piXGal}

 $\pi_1 (X_{Gal}^{Aff})$  is isomorphic to the kernel of $\psi : \pitil \rightarrow S_{18}$.

 \end{theorem}

We denote ${\A} =\pi_1 (X_{Gal}^{Aff}) $.

The exact sequence (\ref{exact1}) gets the form 
\begin{equation}\label{exact2}
1 \rightarrow \A = \ker \psi \rightarrow \pitil  = \frac{\pi_1(\C^2-S,M)}{\sg{\Gamma_j^2, \Gamma_{j'}^2}}
 \stackrel{\psi}{\rightarrow} S_{18} \rightarrow 1.
\end{equation}

The plan is to use the \vK\ Theorem and the braid monodromy
technique to obtain a presentation of  $\pi_1(\C^2 - S, M)$ and by adding the relations $\G_j^2=\G_{j'}^2= 1$ to get a presentation of $\pitil$. 
Then we use the Reidemeister Schreier method to obtain a presentation of $\A$ 
(see [7]).

\subsection{Braid monodromy}\label{braid}
\indent

Recall that $S$ is the branch curve of $X$, deg$S=54$. Recall the above 
$N, M', K(x), u, \C_u, E, D$.  
 Let $B_{54}[D,\C_u]$ be the braid group, and $H_1 , \cdots , H_{53}$ be its frame.  Let $M \in \partial D$ and $\pi_1(\C^2 - S, M)$ is the fundamental group of the complement of  $S$, with a g-base $\G_1, \cdots, \G_{27'}$.

The braid monodromy of $S$ is a map $\varphi: \pi_1(E - N, u) \rightarrow B_{54}[D,\C_u]$ defined as follows:
every loop in $E - N$ starting at $u$ has liftings to a system of ${54}$ paths in $(E - N) \times D$ starting at $q_1, \cdots ,q_{27'}$.
 Projecting them to $D$ we get ${54}$ paths in $D$ defining a motion $\{q_1(t), \cdots , q_{27'}(t)\}$ of ${54}$ points in $D$ starting and ending at $\C_u, 0 \leq t \leq 1$. This motion defines a braid in $B_{54} [D , \C_u]$.

\label{th13}
\noindent
\begin{theorem}{\bf The Artin Theorem}

 Let $S$ be a curve and let $\delta _1, \cdots, \delta _q$ be a g-base 
of $\pi _1 (E -N, u)$.
Assume that the singularities of $S$ are cusps, nodes and branch points. 
Let $\varphi: \pi_1(E - N,u) \rightarrow B_{54}[D, \C_u]$ 
be the braid monodromy.
Then for all $i$, there exist a halftwist $V_i \in B_{54}[D, \C_u]$ 
and $r_i \in \Z$, \st $\varphi  (\dl _i)  = V_i ^{r_i}$
and $r_i$ depends on the type of the singularity:
$r_i =  1,2,3,$ for  branch point, node,  cusp respectively.
\end{theorem}

\begin{proposition}{\bf [17, Proposition VI. 2.1]}\label{pr14}

 Let $S$  be the regenerated branch curve of degree $54$ in $\C \P^2$.  Let $\pi, u, D, E, \C_u$ be as above.  Let $\varphi$ be the braid
 monodromy of $S$ \wrt   $\pi, u.$  Let $\delta_1, \cdots ,
 \delta_q$ be a g-base of $\pi_1(E - N,u)$.  Then
$ \prodl^q_{i=1} \varphi(\delta_i) = \Delta^2_{54}$.
 \end{proposition}

\section{$\Delta^2_{54} =
\prod\limits^9_{i=1} C_i H_{V_i}$}\label{sec:310}

\subsection{$C_i$}
\indent

There are some lines in $X_0$ which do not meet , but when projecting 
them to $\C\P^2$, they may intersect. These intersections $\tC_i$, $i = 1,
\cdots , 9$ are called parasitic intersections. 

Recall from [6, Section 5]: $\tilde{C}_1 = \prod\limits_{t=1,2,4,6,13,22} 
D_t \ , \ \tilde{C}_2 = \prod\limits_{t=3,7,9,14,17} D_t \ , 
\tilde{C}_3 = \prod\limits_{t = 5,10,18,23} D_t \ ,\\
 \ \tilde{C}_4 = \prod\limits_{t=8,11,15,19} D_t \ ,\ \tilde{C}_5 = \prod\limits_{t=12,20,24} D_t \ ,  \tilde{C}_6 = \prod\limits_{t=16,25} D_t \ ,  \ \tilde{C}_7 = \prod\limits_{t=21,26} D_t \ , 
\tilde{C}_8 = D_{27} \ , \ \tilde{C}_9 = Id$.

$D_t ,  1 \leq t \leq 27$,  are regenerated.  We use the
Complex Conjugation [6, Subsection 7.2] and  obtain the following results 
(denoted as above):\\
\\
$D_1 = D_2 = D_3 = Id \ , \ D_4 = Z^2_{33',44'} \ , \
D_5 = \stackrel{\scriptstyle (4)(4')}{\uZ^2}_{\hspace{-.3cm}11',55'} \ ,  \
 \ D_6 = \uZ^2_{33',66'} \cdot \ Z^2_{55',66'}$ ,\\
$D_7 = \prod\limits_{i=2,4} \stackrel{\scriptstyle (6)(6')}{\uZ^2}_{\hspace{-.3cm}ii',77'} \cdot \ \uZ^2_{55',77'} \ ,
\ D_8 = \prod\limits^3_{i=1} \stackrel{\scriptstyle (6)-(7')}{\uZ^2}_{\hspace{-.3cm}ii',88'} \ , \   D_9 = \prod\limits_{i=2, 4-6,8} \uZ^2_{ii',99'} \ , \ \\
D_{10} = \prod\limits_{i=1,4,6,7} \stackrel{\scriptstyle (9)(9')}
{\uZ^2}_{\hspace{-.3cm}ii',10 \; 10'} \cdot \ \uZ^2_{88',10 \; 10'} \ ,
\ D_{11} = \prod\limits_{i=1-3, 6,7} \stackrel{\scriptstyle (9)-(10')}{\uZ^2}_{\hspace{-.4cm}ii',11 \; 11'} \ ,\\
D_{12} = \prod\limits^5_{i=1} \stackrel{\scriptstyle (9)-(11')}
{\uZ^2}_{\hspace{-.4cm}ii',12 \; 12'} \ , \
D_{13} = \prod\limits^{12}_{\stackrel{i=3}{i \neq 4,6}}   \uZ^2_{ii', 13 \; 13'} \ , \
 D_{14} = \prod\limits^{11}_{\stackrel{i=2}{i \neq 3,7,9}}
\stackrel{\scriptstyle (13)(13')}{\uZ^2}_{\hspace{-.4cm}ii',14 \; 14'} \cdot \
\uZ^2_{12 \; 12', 14 \; 14'} \ , \\
D_{15} = \prod\limits^{10}_{\stackrel{i=1}{i \neq 4,5,8}} \stackrel{\scriptstyle (13)-(14')}
{\uZ^2}_{\hspace{-.4cm}ii',15 \; 15'} \cdot \
\uZ^2_{12 \; 12' ,\; 15 \; 15'} \ , \
D_{16} = \prod\limits^{8}_{i=1} \stackrel{\scriptstyle (13)-(15')}
{\uZ^2}_{\hspace{-.4cm}ii',16 \; 16'} \ , \
D_{17} = \prod\limits^{16}_{\stackrel{i=2}{i \neq 3,7,9,14}}
\uZ^2_{ii',17 \; 17'} \ , \ \\
D_{18} = \prod\limits^{15}_{\stackrel{i=1}{i \neq 2,3,5,10}}
\stackrel{\scriptstyle (17)(17')}
{\uZ^2}_{\hspace{-.4cm}ii',18 \; 18'} \cdot \ \uZ^2_{16 \; 16', \; 18 \; 18'}
\ , \
D_{19} = \prod\limits^{14}_{\stackrel{i=1}{i \neq 4,5,8,11}}
\stackrel{\scriptstyle (17)-(18')}
{\uZ^2}_{\hspace{-.4cm}ii',19 \; 19'} \cdot \ \uZ^2_{16 \; 16', \; 19 \; 19'}
 \ , \ \\
D_{20} = \prod\limits^{15}_{\stackrel{i=1}{i \neq 6-8,12}}
\stackrel{\scriptstyle (17)-(19')}
{\uZ^2}_{\hspace{-.4cm}ii',20 \; 20'} \cdot \ \uZ^2_{16 \; 16', \; 20 \; 20'}
 \ , \
D_{21} = \prod\limits^{12}_{i=1}
\stackrel{\scriptstyle (17)-(20')}
{\uZ^2}_{\hspace{-.3cm}ii',21 \; 21'} \ , \
D_{22} = \prod\limits^{21}_{\stackrel{i=3}{i \neq 4,6,13}}
\uZ^2_{ii' , 22 \; 22'}\ , \\$
$D_{23} = \prod\limits^{20}_{\stackrel{i=1}{i \neq 2,3,5,10,18}}
\stackrel{\scriptstyle (22)(22')}
{\uZ^2}_{\hspace{-.4cm}ii',23 \; 23'} \cdot  \uZ^2_{21 \; 21', 23 \; 23'} \ , \
D_{24} = \prod\limits^{19}_{\stackrel{i=1}{i \neq 6-8, 12}}
\stackrel{\scriptstyle (22)-(23')}
{\uZ^2}_{\hspace{-.4cm}ii',24 \; 24'} \cdot  \ \uZ^2_{21 \; 21', 24 \; 24'} \ ,\\
D_{25} = \prod\limits^{20}_{\stackrel{i=1}{i \neq 9-12,16}}
\stackrel{\scriptstyle (22)-(24')}
{\uZ^2}_{\hspace{-.4cm}ii',25 \; 25'} \cdot \ \uZ^2_{21 \; 21', 25 \; 25'} \ , \
D_{26} = \prod\limits^{19}_{\stackrel{i=1}{i \neq 13-16}}
\stackrel{\scriptstyle (22)-(25')}
{\uZ^2}_{\hspace{-.4cm}ii',26 \; 26'} \cdot  \
\stackrel{\scriptstyle (21)(21')}
{\uZ^2}_{\hspace{-.4cm}20 \; 20', 26 \; 26'} \ , \  \\
D_{27} = \prod\limits^{16}_{i=1}
\stackrel{\scriptstyle (22)-(26')}
{\uZ^2}_{\hspace{-.4cm}ii',27 \; 27'} $.

During the regeneration, each $\tilde{C}_i$ is regenerated to $C_i$, 
$1 \leq i \leq 9$.

Each $C_i$ is now a product of the certain regenerated $D_t$ 
(as shown for $\tilde{C}_i$).

\subsection{$H_{V_i}$}
\indent

Recall that the regenerated branch curve compounds nine 6-points ${V_i}$, 
$i=1, \cdots, 9$. We 
regenerate in their neighbourhood. $H_{V_i}$ are the resulting 
regenerated braid monodromies.  
We computed in detail $H_{V_i}$ for $i=1, 4, 7$ in [6, Sections 8, 9, 10].
In this section we present the resulting  
$H_{V_i}$ 
and the paths which correspond to their different factors.
We have proved in [4 ] an Invariance Property to each  $H_{V_i}$, therefore in the following tables , one can see expressions such as $\rho^j_4 \rho^i_3 z_{34} \rho^{-i}_3 \rho^{-j}_4$ where $\rho^j_4=Z_{44'}^j$ and  $\rho^i_4=Z_{33'}^i$,  $i, j \in \Z$. 
 The meaning is that we can consider any relation we need for calculations, but setting $i=j=\pm 1 $ is enough in order to derive all other possibilities.   

There are two tables which correspond to each $H_{V_i}$.The first table 
presents the paths which correspond to the factors in $H_{V_i}$ and some complex conjugates, the braids themselves and the exponent of the braids according to the Emil Artin Theorem.The second table presents the same ones but for the factors of the form $\hF_1 \cdot (\hF_1)^{\rho^{-1}}$.\\

\bigskip
\noindent
$H_{V_1} = (\uZ^2_{4'i})_{i = 22',55',6,6'} \cdot Z^3_{33',4} \cdot 
(\uZ^2_{4i})_{i = 22',55',6,6'}^{Z^2_{33',4}} \cdot
\uZ^3_{11',4} \cdot 
(Z_{44'})^{Z^2_{33',4}\uZ^2_{11',4}} \cdot
Z^3_{55',6} \cdot (\uZ^2_{i6})^{Z^2_{i,22'}Z^2_{33',4}}_{i=11',33'} \cdot
\uZ^3_{22',6} \cdot
(\uZ^2_{i6'})^{\uZ^2_{i6}\uZ^2_{i,55'}Z^2_{i,22'}Z^2_{33',4}}_{i = 11',33'} \cdot
(Z_{66'})^{Z^2_{55',6}\uZ^2_{22',6}} \cdot       
\left(\hF_1(\hF_1)^{Z_{33'}^{-1}Z_{55'}^{-1}}\right)^{Z^2_{55',6}Z^2_{33',4}}$.

\begin{center}
\begin{tabular}[hb]{|l|p{3in}|l|c|}\hline
& The paths/complex conjugates & The braids & The exponent of braids  \\ [-.2cm] \hline 
& & & \\ [-.5cm]
(1) & $\vcenter{\hbox{\epsfbox{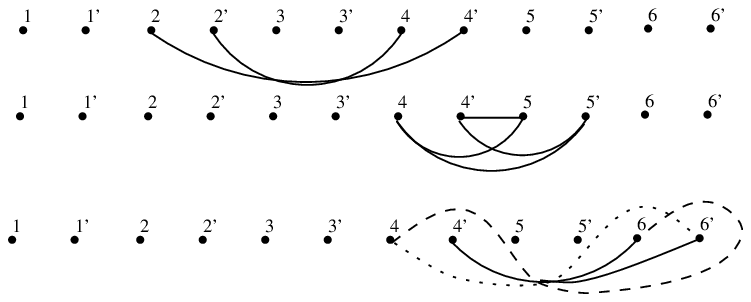}}}$ & $\begin{array}{lll}\rho^j_4 \rho^i_m \underline{z}_{m4} \rho^{-i}_m \rho^{-j}_4 \\ [-.2cm]
m = 2,5  \\  [-.2cm] 
\rho^i_6 \rho^i_4 \underline{z}_{4'6} \rho^{-i}_4 \rho^{-i}_6
\end{array}$ & 2 \\ \hline
& & & \\ [-.5cm]
(2) & $\vcenter{\hbox{\epsfbox{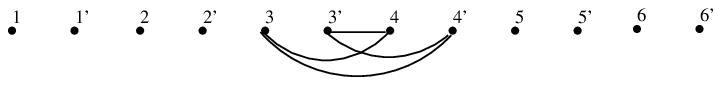}}}$ & $\rho^j_4 \rho^i_3 z_{34} \rho^{-i}_3 \rho^{-j}_4$ & 3  \\[.2cm]  \hline
(3) & $\vcenter{\hbox{\epsfbox{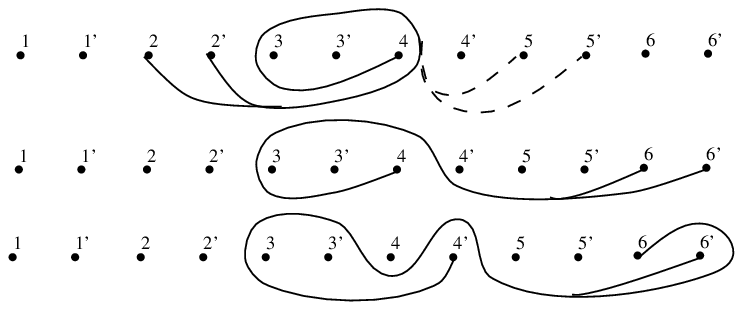}}}$ & $\begin{array}{lll}\rho^j_4 \rho^i_m\tilde{z}_{4m}\rho^{-i}_m \rho^{-j}_4   \\ [-.2cm]
m = 2,5    \\  [-.2cm] 
\rho^i_6 \rho^i_4 \tilde{z}_{46} \rho^{-i}_4 \rho^{-i}_6
\end{array}$ &  2 \\ [1.5cm]\hline
(4) &  $\vcenter{\hbox{\epsfbox{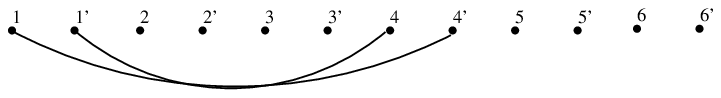}}}$ & $\rho^j_4 \rho^i_1 \underline{z}_{14} \rho^{-i}_1 \rho^{-j}_4$ & 3 \\ [.3cm] \hline
(5) & $\vcenter{\hbox{\epsfbox{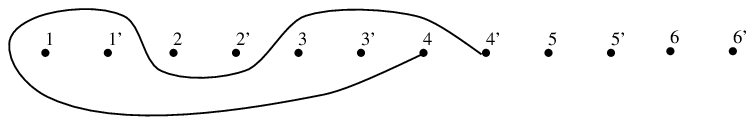}}}$ & $\tilde{z}_{44'}$  & 1 \\
[.3cm] \hline 
& & & \\ [-.4cm]
(6) &   $\vcenter{\hbox{\epsfbox{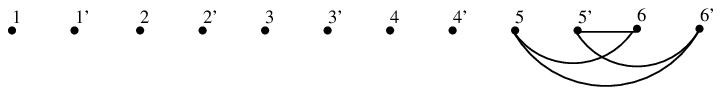}}}$ &  $\rho^j_6 \rho^i_5 z_{56} \rho^{-i}_5 \rho^{-j}_6$ & 3  \\  \hline  
& & & \\  [-.7cm]
(7) & $\vcenter{\hbox{\epsfbox{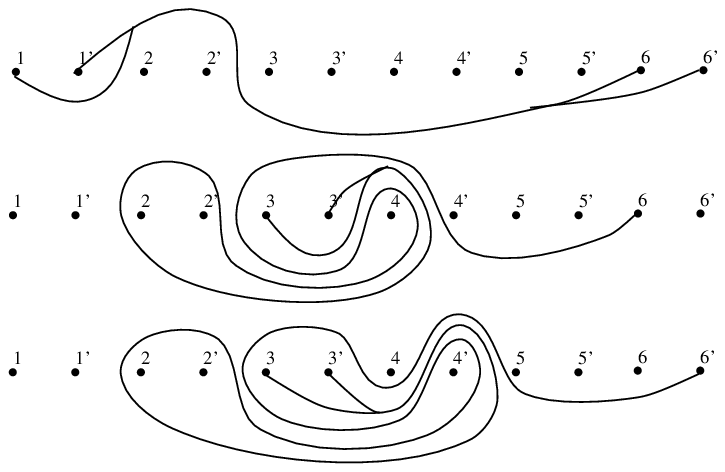}}}$ & $\begin{array}{ll}
\rho^j_6 \rho^i_m \tilde{z}_{m6}\rho^{-i}_m \rho^{-j}_6   \\[-.2cm]
m = 1,3\end{array}$ & 2 \\ \hline
& & & \\ [-.5cm]
(8) &   $\vcenter{\hbox{\epsfbox{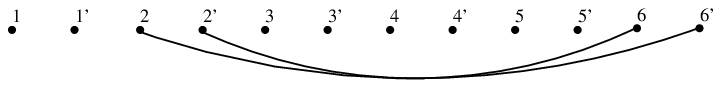}}}$ & $\rho^j_6 \rho^i_2 \underline{ z}_{26} \rho^{-i}_2 \rho^{-j}_6$ & 3 \\  [.3cm]
 \hline
& & & \\ [-.5cm](9) &  $\vcenter{\hbox{\epsfbox{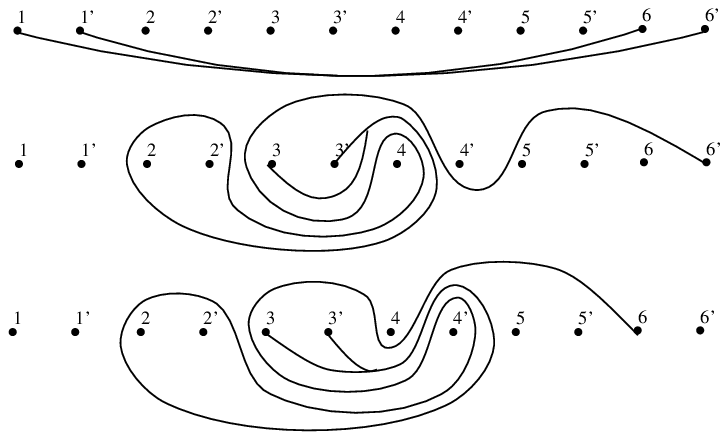}}}$ & $
\begin{array}{ll}
\rho^j_6 \rho^i_m \tilde{z}_{m6}\rho^{-i}_m \rho^{-j}_6\\ [-.2cm]
m = 1,3\end{array}$ & 2 \\ \hline
(10) &   $\vcenter{\hbox{\epsfbox{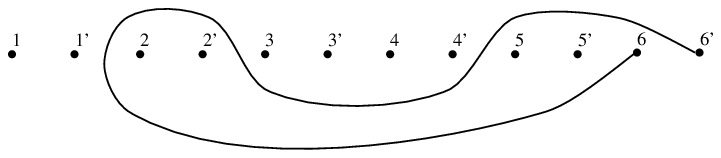}}}$ &  $\tilde{z}_{66'}$ & 1 \\  
   \hline
\end{tabular}
\end{center}
\clearpage
\begin{center} $(\hF_1 \cdot (\hF_1)^{{Z_{33'}^{-1}Z_{55'}^{-1}}})^{Z^2_{55',6}Z^2_{33',4}}$
\end{center}
\begin{center}
\begin{tabular}[H]{|p{3in}|c|c|}\hline
$\vcenter{\hbox{\epsfbox{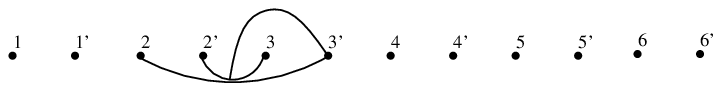}}}$  & $\rho^j_3 \rho^i_2 \underline{z}_{23} \rho^{-i}_2 \rho^{-j}_3$ & 3\\
   \hline
$\vcenter{\hbox{\epsfbox{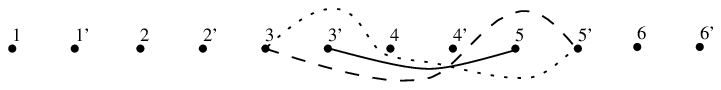}}}$  & $\rho^i_5 \rho^i_3  \underline{z}_{3'5} \rho^{-i}_3 \rho^{-i}_5$ & 2\\
    \hline
$\vcenter{\hbox{\epsfbox{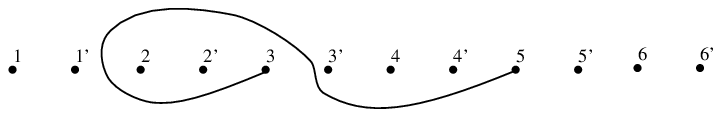}}}$ & $\rho^i_5 \rho^i_3 \tilde{z}_{35} \rho^{-i}_3 \rho^{-i}_5$ & 2\\
$\vcenter{\hbox{\epsfbox{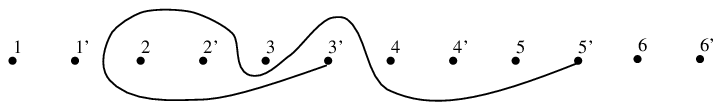}}}$ & & \\
$\vcenter{\hbox{\epsfbox{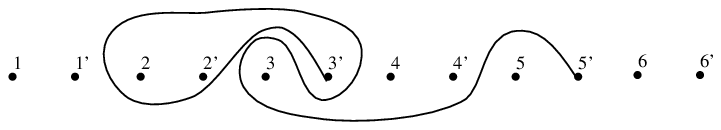}}}$ & & \\ \hline
$\vcenter{\hbox{\epsfbox{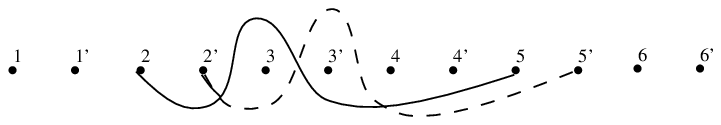}}}$ & $\rho^j_5 \rho^i_2 \tilde{z}_{25} \rho^{-i}_2 \rho^{-j}_5$ & 3\\
$\vcenter{\hbox{\epsfbox{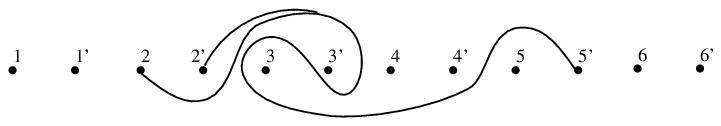}}}$ & & \\ \hline
$\vcenter{\hbox{\epsfbox{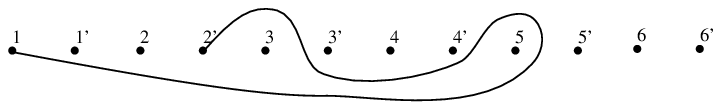}}}$  & $z_{12'}$ & 1\\    \hline
$\vcenter{\hbox{\epsfbox{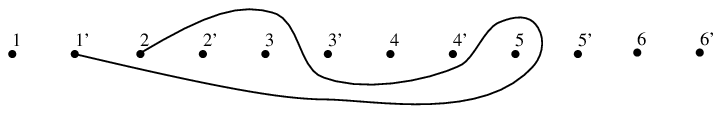}}}$   & $z_{1'2}$ & 1\\
    \hline
$\vcenter{\hbox{\epsfbox{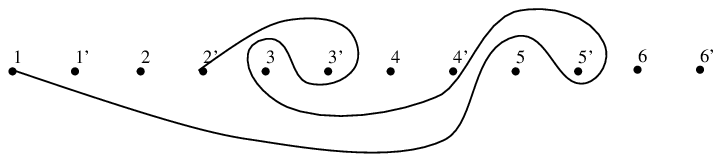}}}$   & $z_{12'}^{\rho^{-1}}$ & 1\\
   \hline
$\vcenter{\hbox{\epsfbox{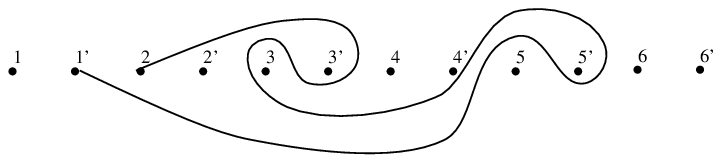}}}$     & $z_{1'2}^{\rho^{-1}}$ & 1
\\   \hline
\end{tabular}
\end{center}
\clearpage
\noindent
$H_{V_2} = Z^3_{33',4} \cdot (\uZ^2_{i4})^{\uZ^{2}_{i,22'}}_{i = 11',5,5',66'}
 \cdot \uZ^3_{22',4} \cdot (\uZ^2_{i4'})^
{\uZ^2_{i4}\uZ^2_{i,33'}\uZ^2_{i,22'}}_{i=11',5,5',66'} \cdot
(Z_{44'})^{Z^2_{33',4}\uZ^2_{22',4}} \cdot Z^3_{5',66'} \cdot (\uZ^2_{33',5})^{Z^2_{33',4}} \cdot \uZ^2_{22',5} \cdot \uZ^3_{11',5} \cdot
(\uZ^2_{33',5'})^{\uZ^2_{33',5}Z^2_{5',66'}Z^2_{33',4}} \cdot
(\uZ^2_{22',5})^{\uZ^{-2}_{11',22'}} \cdot (Z_{55'})^{\uZ^2_{22',5} \uZ^2_{11',5}Z^2_{5',66'}} \cdot \\
\left(\hF_1(\hF_1)^{Z_{33'}^{-1}Z_{66'}^{-1}}\right)^{\uZ^{-2}_{5,66'} Z^{2}_{33',4}}$. \\ 
\begin{center}
\begin{tabular}[hb]{|l|p{3in}|l|c|}\hline
& The paths/complex conjugates & The braids & The exponent of braids \\ [-.2cm]
 \hline
(1) &    $\vcenter{\hbox{\epsfbox{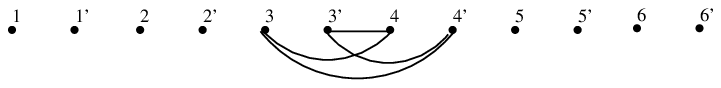}}}$ &  $\rho^j_4 \rho^i_3 \underline{z}_{34} \rho^{-i}_3 \rho^{-j}_4$ & 3\\ 
    \hline
(2) & $\vcenter{\hbox{\epsfbox{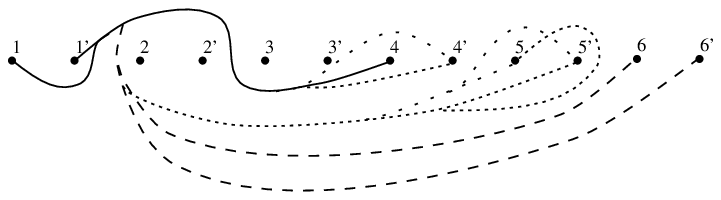}}}$ & $\begin{array}{ll}\rho^j_4 \rho^i_m \tilde{z}_{m4} \rho^{-i}_m \rho^{-j}_4 \\ [-.2cm]
m = 1,6 &  \\ [-.2cm]
\rho^i_5 \rho^i_4 \tilde{z}_{45} \rho^{-i}_4 \rho^{-i}_5 \end{array}$ & 2 \\
  \hline
(3) &  $\vcenter{\hbox{\epsfbox{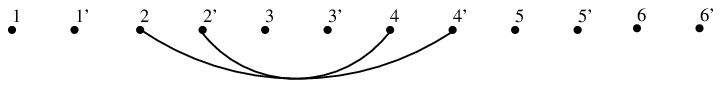}}}$ & $\rho^j_4 \rho^i_2 \underline{z}_{24} \rho^{-i}_2 \rho^{-j}_4$ & 3 \\  
   \hline
(4)  &  $\vcenter{\hbox{\epsfbox{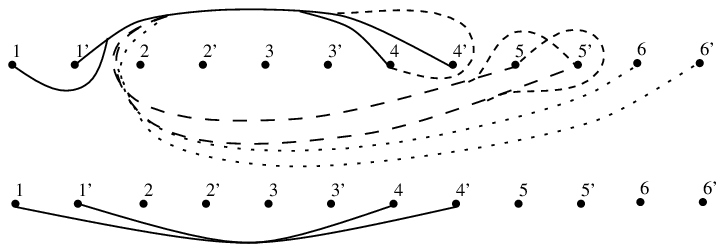}}}$   &   $\begin{array}{ll}
\rho^j_4 \rho^i_m \tilde{z}_{m4} \rho^{-i}_m \rho^{-j}_4\\ [-.2cm]
m = 1,6 &  \\ [-.2cm]
\rho^i_5 \rho^i_4 \tilde{z}_{45} \rho^{-i}_4 \rho^{-i}_5 \end{array}$ & 2 \\
  \hline
(5) &  $\vcenter{\hbox{\epsfbox{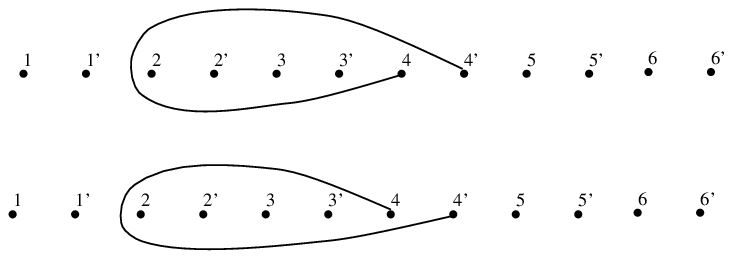}}}$ & $\tilde{z}_{44'}$ & 1\\
    \hline
(6) & $\vcenter{\hbox{\epsfbox{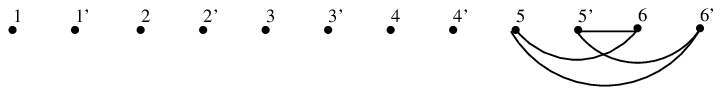}}}$ & $\rho^j_6 \rho^i_5 z_{56} \rho^{-i}_5 \rho^{-j}_6$ & 3  \\ 
   \hline
(7) &  $\vcenter{\hbox{\epsfbox{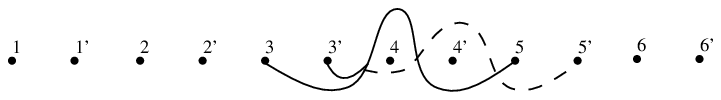}}}$  & $\rho^j_5 \rho^i_3 \tilde{z}_{35}\rho^{-i}_3 \rho^{-j}_5 $ & 2 \\  
  [.5cm] \hline
(8) & $\vcenter{\hbox{\epsfbox{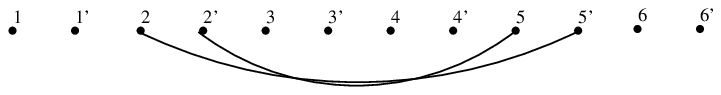}}}$ & $\rho^j_5 \rho^i_2 \underline{z}_{25} \rho^{-i}_2 \rho^{-j}_5$ & 2 \\ 
    \hline
(9) & $\vcenter{\hbox{\epsfbox{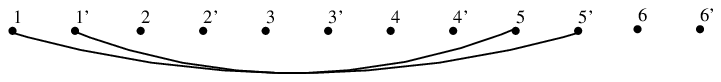}}}$ & $\rho^j_5 \rho_1^i \underline{z}_{15} \rho^{-i}_1 \rho^{-j}_5$ & 3 \\  
   \hline
(10) & $\vcenter{\hbox{\epsfbox{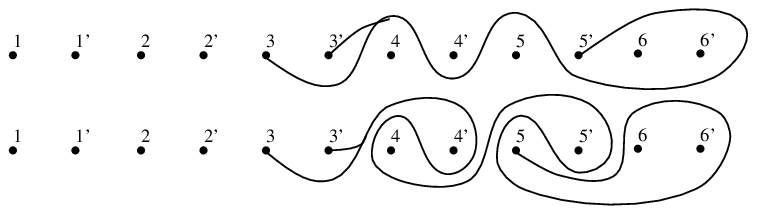}}}$  & $\rho^j_5 \rho^i_3 \tilde{z}_{35'} \rho^{-i}_3 \rho^{-j}_5$ & 2 \\ 
    \hline
(11) &   $\vcenter{\hbox{\epsfbox{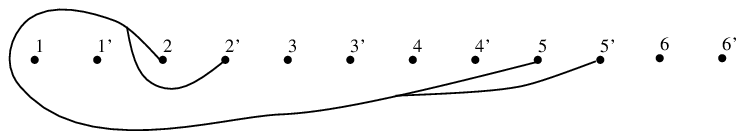}}}$  & $\rho^j_5 \rho^i_2 \tilde{z}_{25} \rho^{-i}_2 \rho^{-j}_5$  & 2 \\  
   \hline
(12) &  $\vcenter{\hbox{\epsfbox{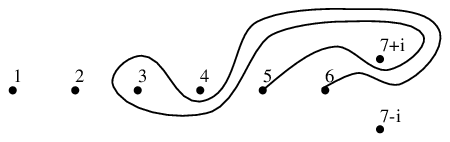}}}$  & $\tilde{z}_{55'}$ & 1 \\ 
  \hline
\end{tabular}
\end{center}

\clearpage
\begin{center}
$(\hF_1 \cdot (\hF_1)^{Z_{33'}^{-1}Z_{66'}^{-1}})^{\uZ^{-2}_{5,66'} Z^{2}_{33',4}}$
\end{center}
\begin{center}
\begin{tabular}[H]{|l|c|c|}\hline
$\vcenter{\hbox{\epsfbox{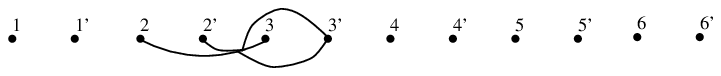}}}$  & $\rho^j_3 \rho^i_2 \underline{z}_{23} \rho^{-i}_2 \rho^{-j}_3$ & 3\\
    \hline
$\vcenter{\hbox{\epsfbox{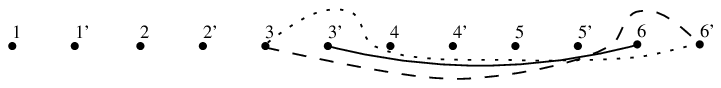}}}$  & $\rho^i_6 \rho^i_3 \underline{z}_{3'6} \rho^{-i}_3 \rho^{-i}_6$ & 2\\
    \hline
$\vcenter{\hbox{\epsfbox{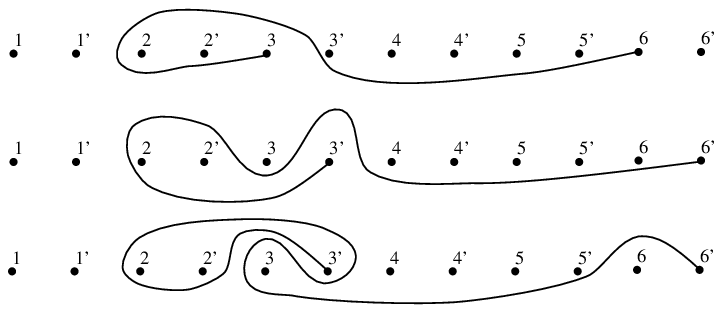}}}$ & $\rho^i_6 \rho^i_3 \tilde{z}_{36} \rho^{-i}_3 \rho^{-i}_6$ & 2\\
   \hline
$\vcenter{\hbox{\epsfbox{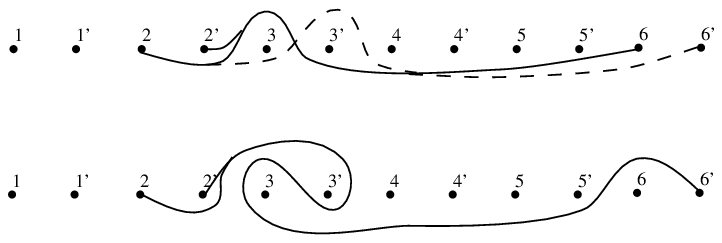}}}$ & $\rho^j_6 \rho^i_2 \tilde{z}_{26} \rho^{-i}_2 \rho^{-j}_6$ & 3\\
    \hline
$\vcenter{\hbox{\epsfbox{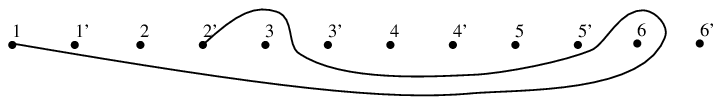}}}$ &  $z_{12'}$ & 1\\
   \hline
$\vcenter{\hbox{\epsfbox{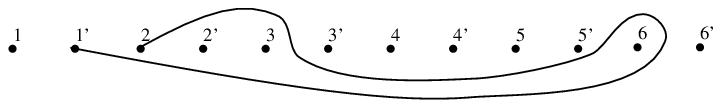}}}$ &  $z_{1'2}$ & 1\\
    \hline
$\vcenter{\hbox{\epsfbox{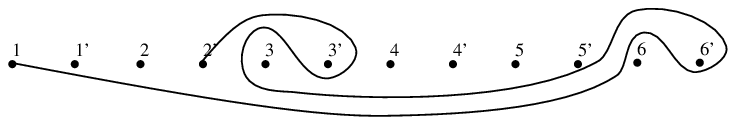}}}$  &  $z_{12'}^{\rho^{-1}}$ & 1\\
    \hline
$\vcenter{\hbox{\epsfbox{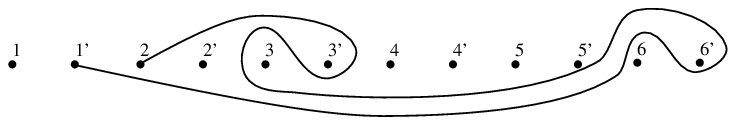}}}$  &  $z_{1'2}^{\rho^{-1}}$ & 1\\
    \hline
\end{tabular}
\end{center}

\clearpage
\noindent
$H_{V_3} = (\uZ^2_{3i})_{i=22',5,5',66'} \cdot 
Z^3_{3',44'} \cdot  
\uZ^3_{11',3} \cdot 
(\uZut_{3'i})^
{\uZ^{2}_{3i}Z^2_{3',44'}}_{i = 22',5,5',66'}
 \cdot 
(Z_{33'})^{\uZ^2_{11',3}Z^2_{3',44'}} \cdot
(\uZ^2_{i5})^{Z^2_{22',3}\uZ^2_{11',3}}_ {i = 11',44'} \cdot
Z^3_{5',66'} \cdot (\uZ^3_{22',5})^{Z^2_{22',3} \uZ^2_{11',3}} \cdot
(\uZ^2_{i5})^{\uZ^{-2}_{22',i} Z^2_{22',3}\uZ^2_{11',3}}_{i = 11',44'} \cdot
(Z_{55'})^{Z^2_{5',66'}Z^2_{44',5}\uZ^2_{22',5}Z^2_{22',3}\uZ^2_{11',3}} \cdot
\left(\hF_1(\hF_1)^{Z_{44'}^{-1}Z_{66'}^{-1}}\right)^{Z^{2}_{22',3} \uZ^{2}_{11',3} 
\uZ^{-2}_{5,66'}}$.\\
\begin{center}
\begin{tabular}[hb]{|l|p{3in}|l|c|}\hline
& The paths/complex conjugates & The braids & The exponent of braids \\ [-.2cm]
  \hline 
(1) &  $\vcenter{\hbox{\epsfbox{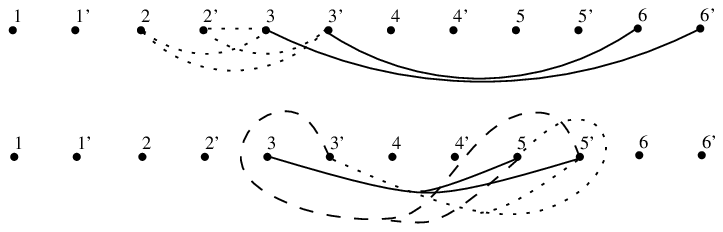}}}$ &  $\begin{array}{lll}
  \rho^j_m \rho^i_3 \underline{z}_{3m} \rho^{-i}_3 \rho^{-j}_m  \\ [-.2cm]
m = 2,6  \\ [-.2cm]
\rho^i_5 \rho^i_3 \underline{z}_{35} \rho^{-i}_3 \rho^{-i}_5 
\end{array}$ & 2 \\
    \hline
(2) & $\vcenter{\hbox{\epsfbox{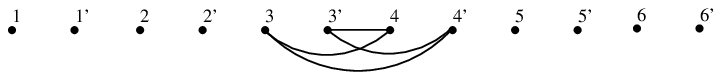}}}$ & $\rho^j_4 \rho^i_3 
\underline{ z}_{3'4} \rho^{-i}_3 \rho^{-j}_4$ & 3  \\   \hline 
&&&\\ [-.5cm]
(3) &  $\vcenter{\hbox{\epsfbox{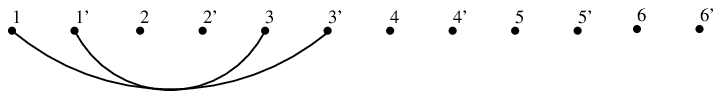}}}$ & $\rho^j_3 \rho^i_1 \underline{z}_{13} \rho^{-i}_1 \rho^{-j}_3$ & 3 \\   
   \hline
(4) &  $\vcenter{\hbox{\epsfbox{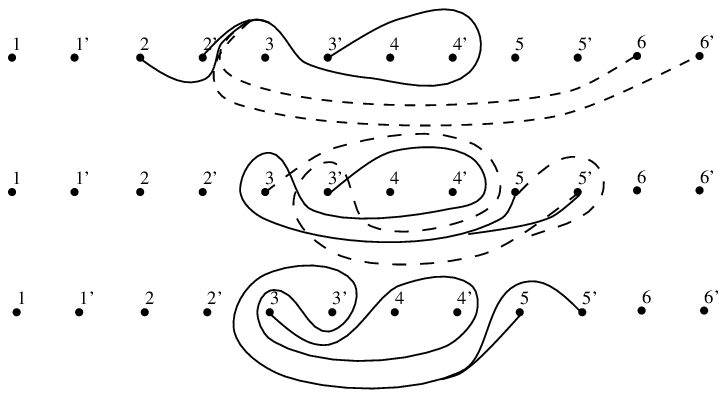}}}$ & $\begin{array}{lll}
\rho^j_m \rho^i_3 z_{3'm}\rho^{-i}_3 \rho^{-j}_m \\ [-.2cm]
m=2,6 \\ [-.2cm]
\rho^i_5 \rho^i_3 \tilde{z}_{3'5} \rho^{-i}_3 \rho^{-i}_5 \end{array}$ & 2 \\
    \hline
(5) &  $\vcenter{\hbox{\epsfbox{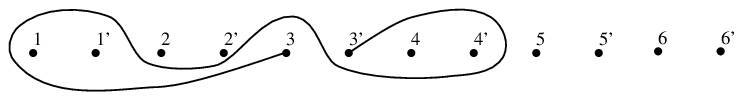}}}$ & $\tilde{z}_{33'} $ & 1\\
   \hline
(6) &  $\vcenter{\hbox{\epsfbox{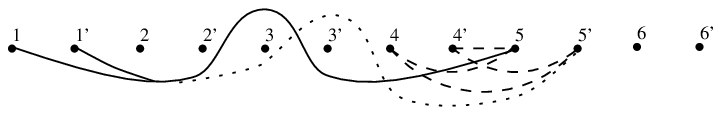}}}$ & $\rho^j_5\rho^i_m \tilde{z}_{m5}\rho^{-i}_m\rho^{-j}_5$ & 2  \\ [-.4cm]
& & $m = 1,4$ &\\ 
   \hline
(7) &  $\vcenter{\hbox{\epsfbox{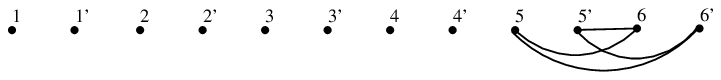}}}$ & $\rho^j_6 \rho^i_5 
\underline{ z}_{5'6}\rho^{-i}_5 \rho^{-j}_6 $ & 3 \\  
   \hline
(8) &  $\vcenter{\hbox{\epsfbox{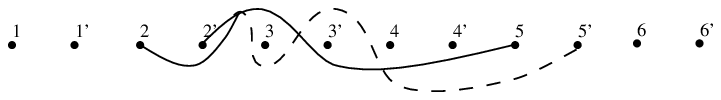}}}$ & $\rho^j_5 \rho^i_2 \tilde{z}_{25} \rho^{-i}_2 \rho^{-j}_5$ & 3 \\ 
    \hline
(9) &  $\vcenter{\hbox{\epsfbox{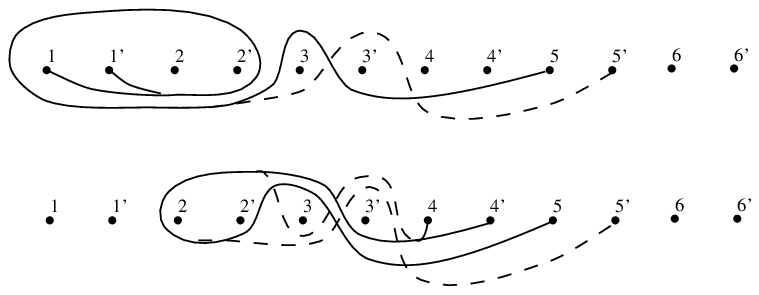}}}$ & $\begin{array}{lll}\rho_5^j \rho_m^i \tilde{z}_{m5} \rho_m^{-i} \rho_5^{-j} \\ [-.2cm]
m = 1,4\end{array}$ & 2\\
   \hline
(10) & $\vcenter{\hbox{\epsfbox{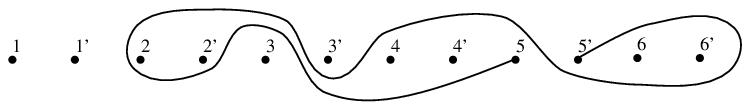}}}$ & $\tilde{z}_{55'}$ & 1 \\  
   \hline
\end{tabular}
\end{center}
\clearpage
\begin{center}
$(\hF_1 \cdot (\hF_1)^{Z_{44'}^{-1}Z_{66'}^{-1}})^{Z^{2}_{22',3} \uZ^{2}_{11',3} 
\uZ^{-2}_{5,66'}}$ 
\end{center}
\begin{center}
\begin{tabular}[H]{|l|c|c|}\hline
$\vcenter{\hbox{\epsfbox{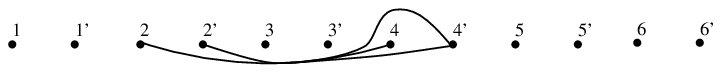}}}$ &  $\rho^j_4 \rho^i_2 \underline{z}_{24} \rho^{-i}_2 \rho^{-j}_4$ & 3\\
    \hline
$\vcenter{\hbox{\epsfbox{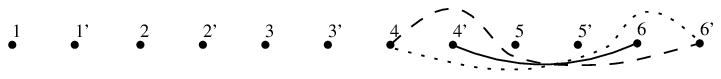}}}$  &   $\rho^i_6 \rho^i_4 \underline{z}_{4'6} \rho^{-i}_4 \rho^{-i}_6$ & 2\\
  \hline
$\vcenter{\hbox{\epsfbox{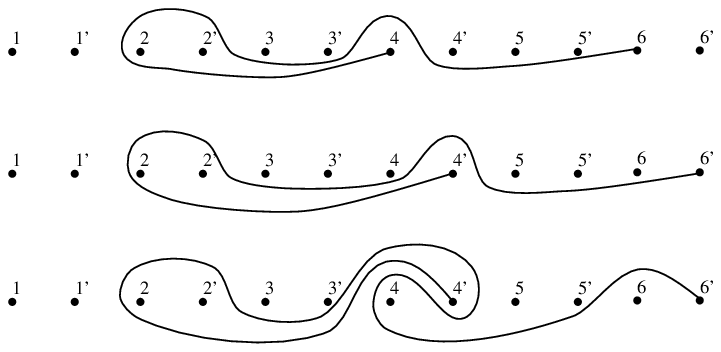}}}$  & $\rho^i_6 \rho^i_4 \tilde{z}_{46} \rho^{-i}_4 \rho^{-i}_6$ & 2\\
    \hline
$\vcenter{\hbox{\epsfbox{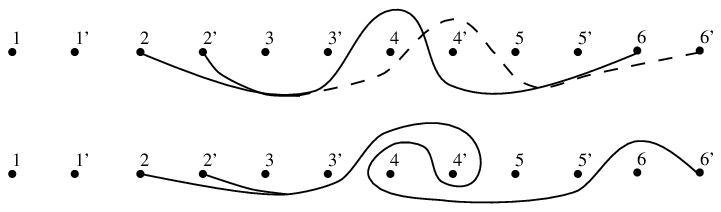}}}$ & $\rho^j_6 \rho^i_2 \tilde{z}_{26} \rho^{-i}_2 \rho^{-j}_6$ & 3\\
    \hline
$\vcenter{\hbox{\epsfbox{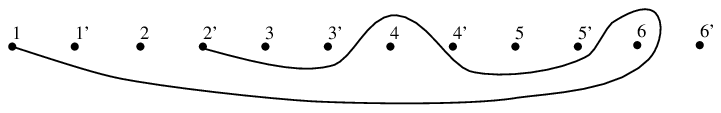}}}$    &  $z_{12'}$ & 1\\
   \hline
$\vcenter{\hbox{\epsfbox{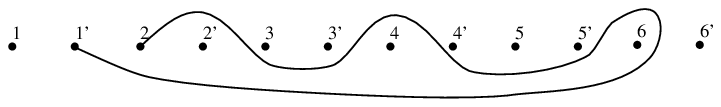}}}$  &  $z_{1'2}$ & 1\\
    \hline
$\vcenter{\hbox{\epsfbox{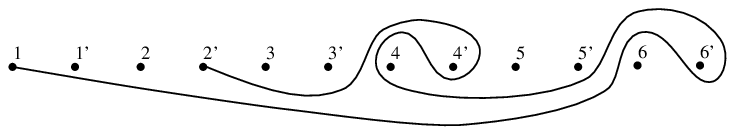}}}$  &  $z_{12'}^{\rho^{-1}}$ & 1\\
    \hline
$\vcenter{\hbox{\epsfbox{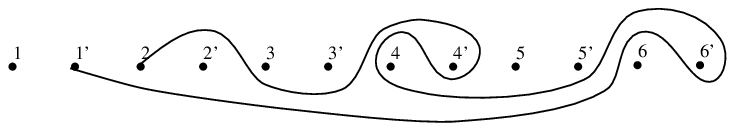}}}$  &  $z_{1'2}^{\rho^{-1}}$ & 1\\
    \hline
\end{tabular}
\end{center}

\clearpage
\noindent
$H_{V_4}$ = $(\uZ^2_{2'i})_{i=33',55',6,6'} \cdot
Z^3_{11',2} \cdot
\bZ^3_{2',44'}
\cdot
(\uZ^2_{2'i})^{\uZ^{-2}_{i,44'}}_{i=33',55',6,6'}
\cdot
(Z_{22'})^{Z^2_{11',2}\uZut_{2',44'}Z^2_{2',33'}}\cdot
(Z_{66'})^{\uZumt_{33',6}  Z^{-2}_{55',6}}\cdot
(\uZ^2_{i6'})^{Z^2_{11',2} }_{i = 11',44'}
\cdot
(\uZ^3_{33',6})^{Z^{-2}_{55',6}}
\cdot
(\uZ^2_{i6})^{Z^{-2}_{55',6}Z^2_{11',2}}_{i = 11',44'}
\cdot Z^3_{55',6} 
\cdot (\hF_1(\hF_1)^{Z_{11'}^{-1}Z_{55'}^{-1}})^{Z^2_{11',2}Z^{-2}_{55',6}}$. \\ 
\begin{center}
\begin{tabular}[hb]{|l|p{3in}|l|c|}\hline
& The paths/complex conjugates & The braids & The exponent of braids \\ [-.2cm] \hline
(1) &  $\vcenter{\hbox{\epsfbox{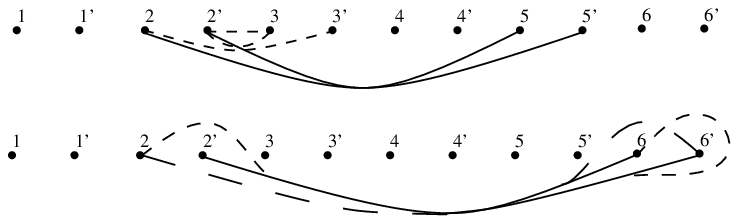}}}$ 
 & $\begin{array}{lll}\rho^j_m \rho^i_2 \underline{z}_{2'm} \rho^{-i}_2 \rho^{-j}_m\\ [-.2cm]
m = 3,5 \\ [-.2cm]
\rho^i_6 \rho^i_2 \underline{z}_{2'6} \rho^{-i}_2 \rho^{-i}_6 
\end{array}$ & 2 \\
    \hline
(2) &  $\vcenter{\hbox{\epsfbox{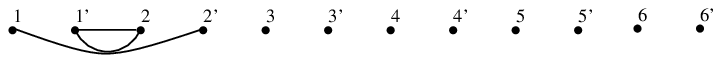}}}$ & $\rho^j_2 \rho^i_1 \underline{z}_{12} \rho^{-i}_1 \rho^{-j}_2$ & 3  \\ 
   \hline
(3) &  $\vcenter{\hbox{\epsfbox{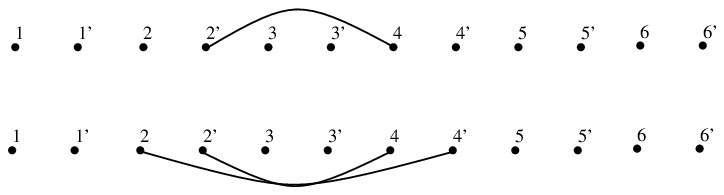}}}$ & $ \rho^j_4 \rho^i_2 \underline{z}_{2'4} \rho^{-i}_2 \rho^{-j}_4$ & 3  
\\  
   \hline
(4) &  $\vcenter{\hbox{\epsfbox{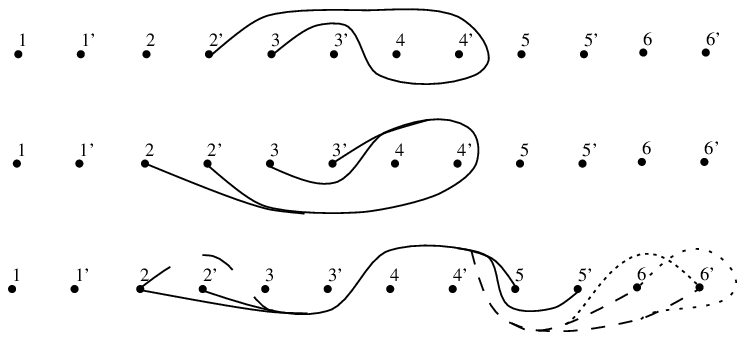}}}$ & $\begin{array}{lll}
\rho^j_m \rho^i_2 \tilde{z}_{2'm} \rho^{-i}_2 \rho^{-j}_m \\ [-.2cm]
m = 3,5\\ [-.2cm]
\rho^i_6 \rho^i_2 \tilde{z}_{2'6} \rho^{-i}_2 \rho^{-i}_6 \end{array}$ & 2 \\
    \hline
(5) & $\vcenter{\hbox{\epsfbox{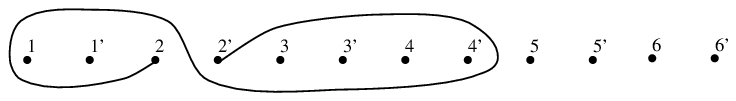}}}$ & $\tilde{z}_{22'}$ & 1\\
   \hline
(6) & $\vcenter{\hbox{\epsfbox{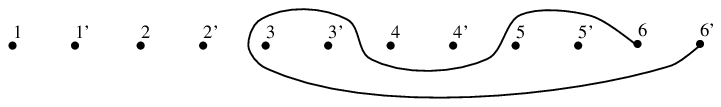}}}$ & $\tilde{z}_{66'}$ & 1  \\ 
   \hline
(7) & $\vcenter{\hbox{\epsfbox{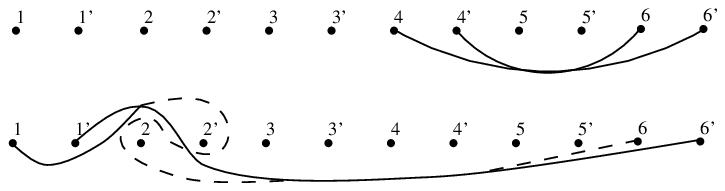}}}$ & $\begin{array}{ll}\rho^j_6 \rho^i_m \tilde{z}_{m6'}\rho^{-i}_m \rho^{-j}_6 \\ [-.2cm]
m = 1,4\end{array}$ & 2\\
   \hline
(8) & $\vcenter{\hbox{\epsfbox{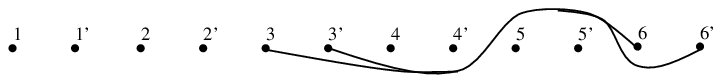}}}$ & $\rho^j_6 \rho^i_3 \tilde{z}_{36} \rho^{-i}_3 \rho^{-j}_6$ & 3 \\ 
    \hline
(9) & $\vcenter{\hbox{\epsfbox{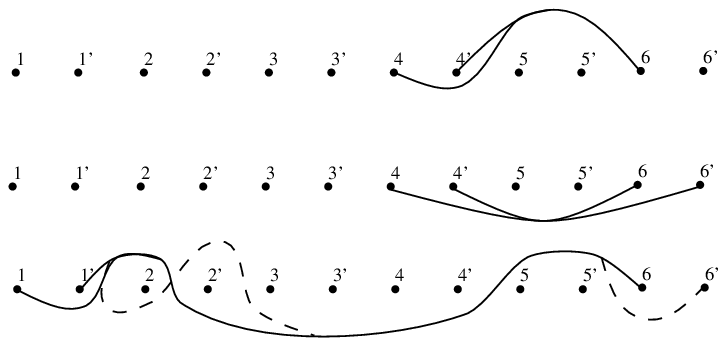}}}$ & 
$\begin{array}{ll}\rho^j_6 \rho^i_m \tilde{z}_{m6} \rho^{-i}_m \rho^{-j}_6
\\ [-.2cm]
m = 1,4 \end{array}$ & 2 \\
   \hline
(10) & $\vcenter{\hbox{\epsfbox{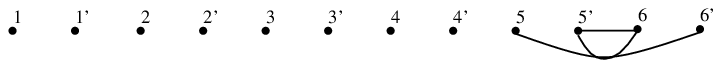}}}$ & 
$\rho^j_6 \rho^i_5 z_{56} \rho^{-i}_5 \rho^{-j}_6$ & 3 \\ 
   \hline
\end{tabular}
\end{center}
\clearpage
\begin{center}
$(\hF_1 \cdot (\hF_1)^{Z_{11'}^{-1}Z_{55'}^{-1}})^{Z^2_{11',2}Z^{-2}_{55',6}}$ 
\end{center}
\begin{center}
\begin{tabular}[H]{|l|c|c|} \hline
 $\vcenter{\hbox{\epsfbox{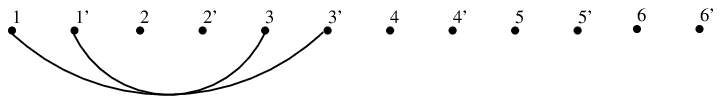}}}$ & $\rho^j_3 \rho^i_1 \underline{z}_{1'3} \rho^{-i}_1 \rho^{-j}_3$ & 3\\ [.3cm]
    \hline
 $\vcenter{\hbox{\epsfbox{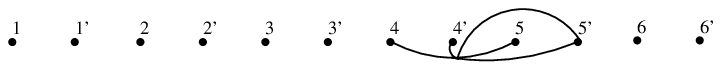}}}$  & $\rho^j_5 \rho^i_4  \underline{z}_{45} \rho^{-i}_4 \rho^{-j}_5$ & 3\\
    \hline
 $\vcenter{\hbox{\epsfbox{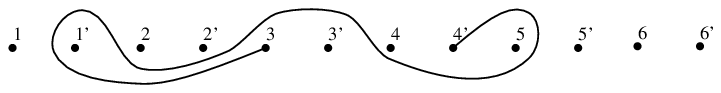}}}$ & $z_{34'}$ & 1\\
    \hline
 $\vcenter{\hbox{\epsfbox{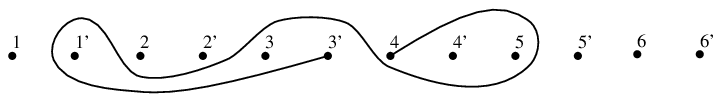}}}$  & $z_{3'4}$ & 1\\
    \hline
 $\vcenter{\hbox{\epsfbox{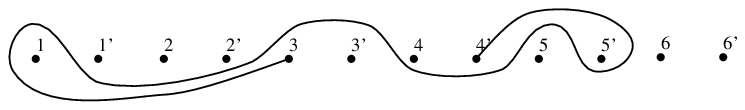}}}$  & $z_{34'}^{\rho^{-1}}$ & 1\\
    \hline
 $\vcenter{\hbox{\epsfbox{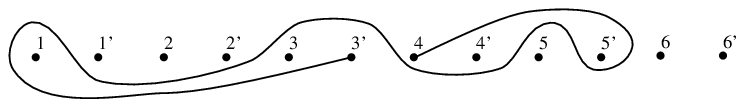}}}$   & $z_{3'4}^{\rho^{-1}}$ & 1\\   \hline
 $\vcenter{\hbox{\epsfbox{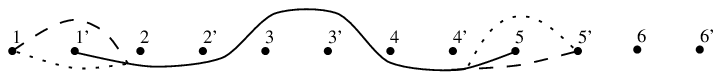}}}$ & $\rho^i_5 \rho^i_1 \tilde{z}_{1'5} \rho^{-i}_1 \rho^{-i}_5$ & 2\\
    \hline
 $\vcenter{\hbox{\epsfbox{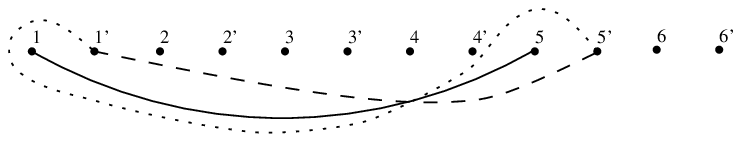}}}$  & $\rho^i_5 \rho^i_1 \underline{z}_{15} \rho^{-i}_1 \rho^{-i}_5$ & 2\\
   \hline
\end{tabular}
\end{center}
\clearpage
\noindent
$H_{V_5} = Z^3_{1',22'} \cdot (Z_{66'})^{\uZumt_{44',6} Z^{-2}_{55',6}} \cdot 
\uZ^2_{33',6'} \cdot \bZ^3_{44',6}\cdot(\uZ^2_{33',6})^{Z^{-2}_{55',6}} \cdot
(\uZ^2_{22',6'})^{Z^2_{1',22'}} \cdot 
(\uZ^2_{22',6})^{Z^2_{1',22'}Z^{-2}_{55',6}}  \cdot \\
 (\hF_1(\hF_1)^{Z_{22'}^{-1}Z_{55'}^{-1}})^{Z^2_{1',22'}Z^{-2}_{55',6}}
\cdot
\uZ^3_{55',6} \cdot (\uZ^2_{1'i})^{Z^2_{1',22'}}_{i = 44',55',6,6'} \cdot 
\bZ^3_{1',33'} \cdot (\uZ^2_{1i})_{i = 44',55',6,6'} \cdot
(Z_{11'})^{\uZ^2_{1',33'}Z^2_{1',22'}}$. \\
\begin{center}
\begin{tabular}[hb]{|l|p{2.9in}|l|c|}\hline
& The paths/complex conjugates & The braids & The exponent of braids \\ [-.2cm]\hline
(1) & $\vcenter{\hbox{\epsfbox{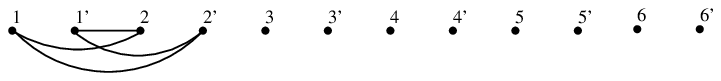}}} $ & $\rho^j_2 \rho^i_1 z_{1'2} \rho^{-i}_1 \rho^{-j}_2$ & 3\\
   \hline
(2) & $\vcenter{\hbox{\epsfbox{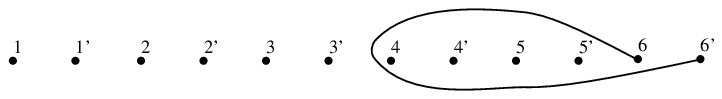}}} $ & $\tilde{z}_{66'}$ & 1  \\ 
   \hline
(3) & $\vcenter{\hbox{\epsfbox{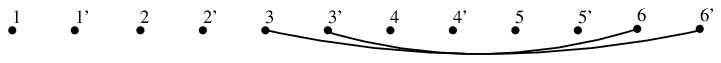}}} $ & $ \rho^j_6 \rho^i_3 \underline{z}_{36'} \rho^{-i}_3 \rho^{-j}_6$ & 2  
\\  
(4) & $\vcenter{\hbox{\epsfbox{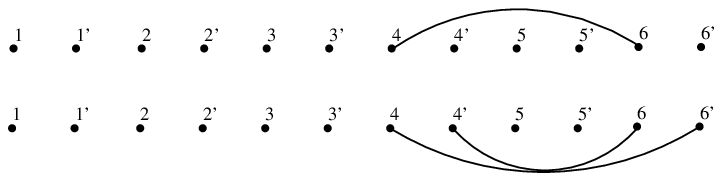}}} $ & $\rho^j_6 \rho^i_4 \underline{z}_{46} \rho^{-i}_4 \rho^{-j}_6$ & 3 \\ 
    \hline
(5) & $\vcenter{\hbox{\epsfbox{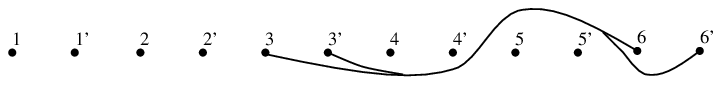}}} $ & $ \rho^j_6 \rho^i_3 \tilde{z}_{36} \rho^{-i}_3 \rho^{-j}_6 $ & 2\\
    \hline
(6) &  $\vcenter{\hbox{\epsfbox{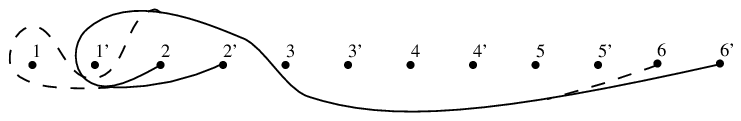}}} $ & $\rho^j_6 \rho^i_2 \tilde{z}_{26'} \rho^{-i}_2 \rho^{-j}_6$ & 2  \\ 
  \hline
(7) &  $\vcenter{\hbox{\epsfbox{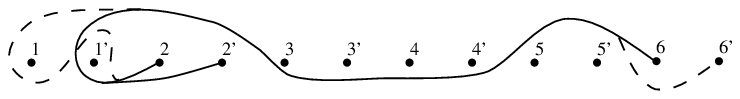}}} $ & $\rho^j_6 \rho^i_2 \tilde{z}_{26}\rho^{-i}_2 \rho^{-j}_6 $ & 2 \\  
   \hline
(8) & $\vcenter{\hbox{\epsfbox{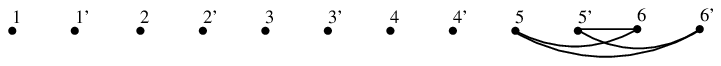}}} $ & $\rho^j_6 \rho^i_5 z_{56} \rho^{-i}_5 \rho^{-j}_6$ & 3 \\
    \hline
(9) & $\vcenter{\hbox{\epsfbox{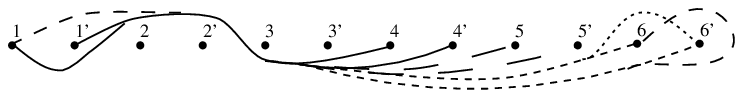}}} $ & $ \begin{array}{lll}\rho^j_m \rho^i_1 \tilde{z}_{1'm} \rho^{-i}_1 \rho^{-j}_m & \\ [-.2cm]
m=4,5 &  \\ [-.2cm]
 \rho^i_6 \rho^i_1 \tilde{z}_{1'6} \rho^{-i}_1 \rho^{-i}_6 \end{array}$ & 2 \\
\hline
(10) & $\vcenter{\hbox{\epsfbox{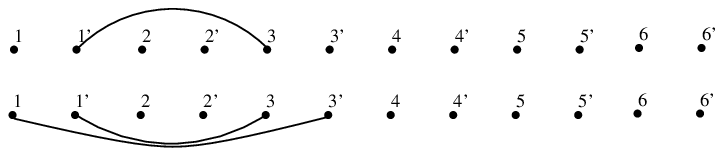}}} $ & $\rho^j_3 \rho^i_1 \underline{z}_{1'3} \rho^{-i}_1 \rho^{-j}_3$ & 3 
\\    \hline
(11) & $\vcenter{\hbox{\epsfbox{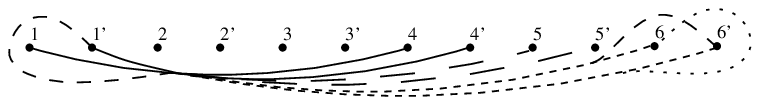}}} $ & $\begin{array}{lll} \rho^j_m \rho^i_1 \underline{z}_{1m} \rho^{-i}_1 \rho^{-j}_m &  \\ [-.2cm]
m=4,5 &   \\ [-.2cm]
 \rho^i_6 \rho^i_1 \underline{z}_{16} \rho^{-i}_1 \rho^{-i}_6 \end{array}$ & 2 \\
\hline
(12) &  $\vcenter{\hbox{\epsfbox{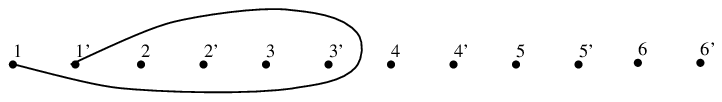}}} $ &  $\tilde{z}_{11'}$  &  1 \\ 
 \hline
\end{tabular}
\end{center}
\clearpage
\begin{center}
$(\hF_1 \cdot (\hF_1)^{Z_{22'}^{-1}Z_{55'}^{-1}})^{\uZ^2_{1',33'}Z^2_{1',22'}}$ 
\end{center}
\vspace{-.1cm}
\begin{center}
\begin{tabular}[H]{|l|c|c|}\hline
$\vcenter{\hbox{\epsfbox{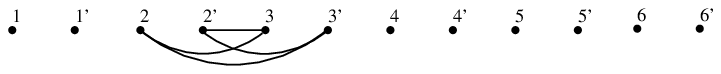}}} $ & $\rho^j_3 \rho^i_2 \underline{z}_{23} \rho^{-i}_2 \rho^{-j}_3$ & 3\\
   \hline
$\vcenter{\hbox{\epsfbox{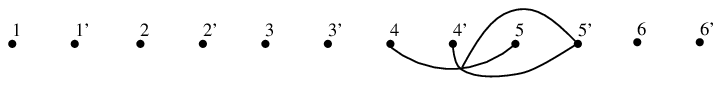}}} $  & $\rho^j_5 \rho^i_4  \underline{z}_{45} \rho^{-i}_4 \rho^{-j}_5$ & 3\\
    \hline
 $\vcenter{\hbox{\epsfbox{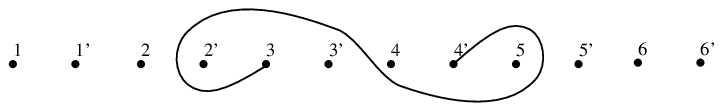}}} $ & $z_{34'}$ & 1\\
       \hline
 $\vcenter{\hbox{\epsfbox{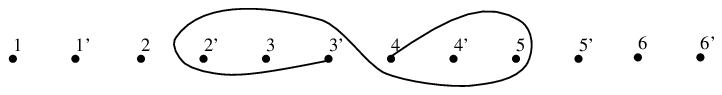}}} $  & $z_{3'4}$ & 1\\
    \hline
 $\vcenter{\hbox{\epsfbox{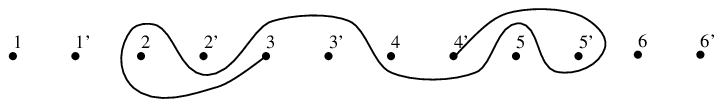}}} $  & $z_{34'}^{\rho^{-1}}$ & 1\\
    \hline
 $\vcenter{\hbox{\epsfbox{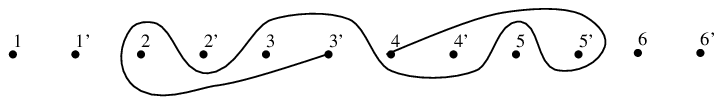}}} $  & $z_{3'4}^{\rho^{-1}}$ & 1 \\ \hline
 $\vcenter{\hbox{\epsfbox{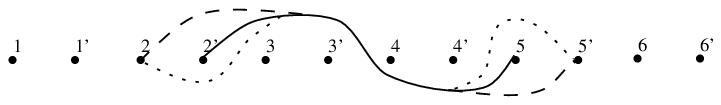}}} $ & $\rho^i_5 \rho^i_2 \tilde{z}_{2'5} \rho^{-i}_2 \rho^{-i}_5$ & 2\\
    \hline
 $\vcenter{\hbox{\epsfbox{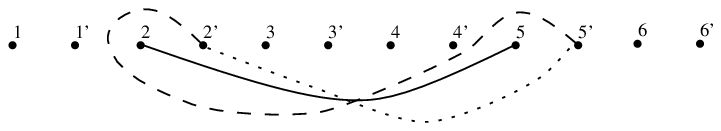}}} $  & $\rho^i_5 \rho^i_2 \underline{z}_{25} \rho^{-i}_2 \rho^{-i}_5$ & 2\\
    \hline
\end{tabular}
\end{center}
\clearpage
\noindent
$H_{V_6} = Z^3_{1',22'} 
 \cdot 
(\uZ^2_{1'i})^{Z^2_{1',22'}}_{i = 33',5,5',66'} \cdot
(\uZ^3_{1',44'})^{Z^2_{1',22'}} \cdot 
(\uZ^2_{1i})_{i = 33',5,5',66'}
\cdot
(Z_{11'})^{\uZ^2_{1',44'}Z^2_{1',22'}}
\cdot (Z_{55'})^{\bZ^{-2}_{33',5}  Z^{-2}_{5',66'}} \cdot
(\uZ^2_{i5'})^{Z^{-2}_{5',66'}Z^2_{1',22'}}_{i = 22',44'} \cdot
\bZ^3_{33',5} \cdot (\uZ^2_{i5})^{Z^{2}_{1',22'}}_{i = 22',44'}  \cdot 
Z^3_{5',66'}\cdot (\hF_1(\hF_1)^{Z_{22'}^{-1}Z_{66'}^{-1}})^{Z^2_{1',22'}Z^{-2}_{5',66'}}$. \\
\vspace{-.68cm}
\begin{center}
\begin{tabular}[hb]{|l|p{3in}|l|c|}\hline
& The paths/complex conjugates & The braids & The exponent of braids \\ [-.2cm]\hline
(1) &  $\vcenter{\hbox{\epsfbox{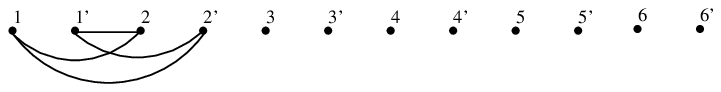}}}$ & $\rho^j_2 \rho^i_1 z_{1'2} \rho^{-i}_1 \rho^{-j}_2$ & 3\\
    \hline
(2) &  $\vcenter{\hbox{\epsfbox{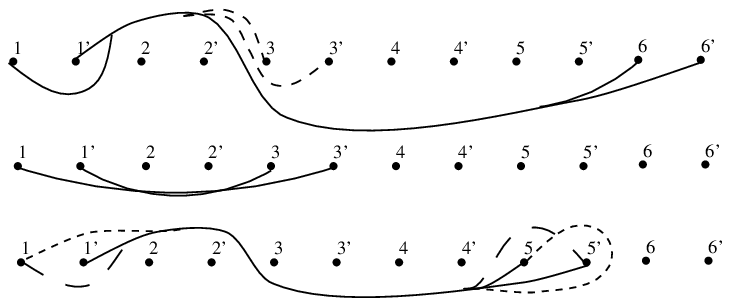}}}$ & $\begin{array}{lll}
\rho^j_m \rho^i_1 \tilde{z}_{1'm} \rho^{-i}_1 \rho^{-j}_m \\ [-.2cm]
m = 3,6 \\ [-.2cm]
\rho^i_5 \rho^i_1 \tilde{z}_{1'5} \rho^{-i}_1 \rho^{-i}_5\end{array}$ & 2 \\
   \hline
(3) & $\vcenter{\hbox{\epsfbox{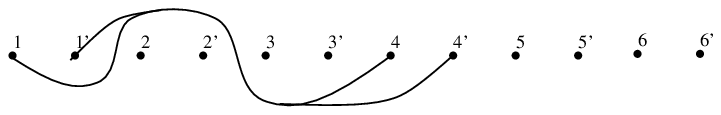}}}$ & $ \rho^j_4 \rho^i_1 \tilde{z}_{1'4} \rho^{-i}_1 \rho^{-j}_4$ & 3  
\\  
   \hline
&&&\\
(4) & $\vcenter{\hbox{\epsfbox{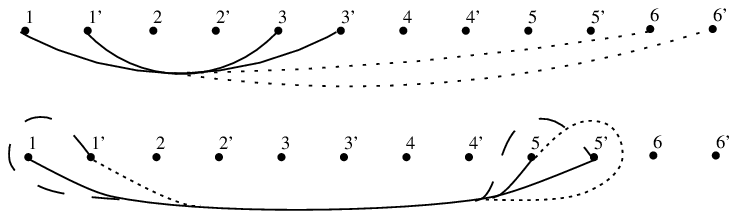}}}$ & $\begin{array}{lll}\rho^j_m \rho^i_1 \underline{z}_{1m} \rho^{-i}_1 \rho^{-j}_m\\ [-.2cm]
 m = 3,6 \\   [-.2cm]
\rho^i_5 \rho^i_1 \underline{z}_{15} \rho^{-i}_1 \rho^{-i}_5\end{array}$ & 2 \\
   \hline
  
(5) & $\vcenter{\hbox{\epsfbox{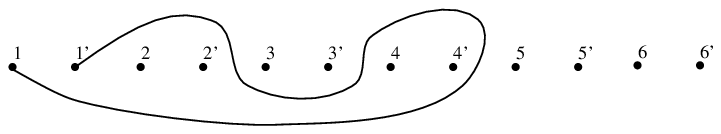}}}$ & $\tilde{z}_{11'} $ & 1\\
    \hline

(6) & $\vcenter{\hbox{\epsfbox{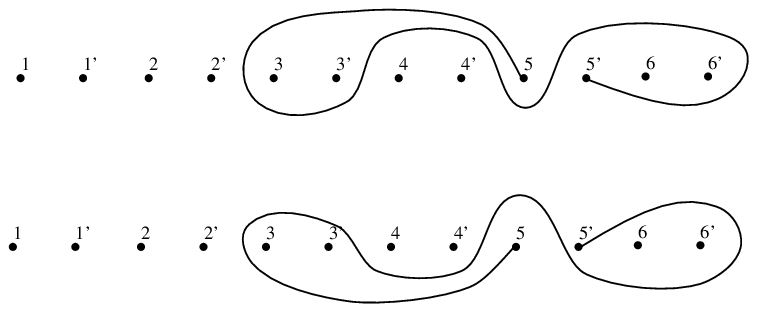}}}$ & $\tilde{z}_{55'}$ & 1  \\

   \hline

(7) & $\vcenter{\hbox{\epsfbox{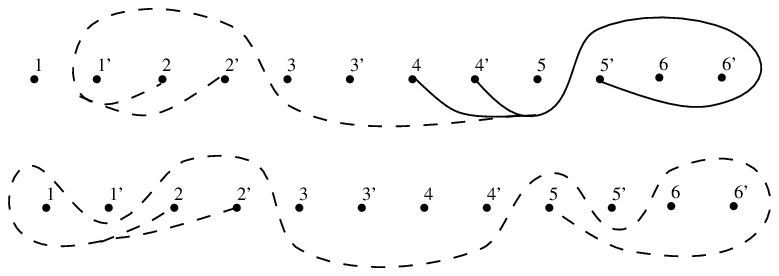}}}$ & $\rho^j_5 \rho^i_m \tilde{z}_{m5'}\rho^{-i}_m \rho^{-j}_5 $ & 2 \\  [-.4cm]
 & & $m = 2,4$ &   \\ \hline

(8) & $\vcenter{\hbox{\epsfbox{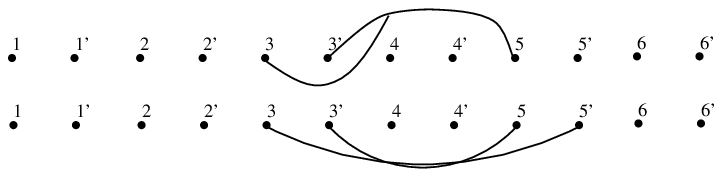}}}$ & $\rho^j_5 \rho^i_3
\underline{z}_{35} \rho^{-i}_3 \rho^{-j}_5$  & 3 \\

    \hline

(9) & $\vcenter{\hbox{\epsfbox{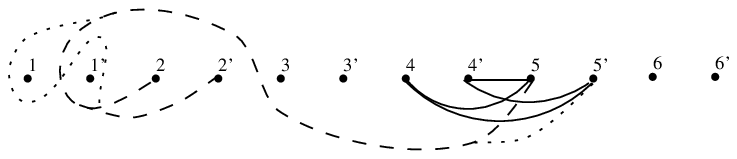}}}$ & $ \rho^j_5 \rho^i_m \tilde{z}_{m5} \rho^{-i}_m \rho^{-j}_5$ & 2 \\ [-.4cm]
   & & $m = 2,4$ &  \\ \hline

(10) & $\vcenter{\hbox{\epsfbox{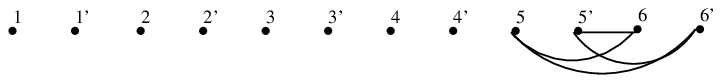}}}$ & $ \rho^j_6 \rho^i_5 z_{5'6} \rho^{-i}_5 \rho^{-j}_6$ & 3 \\

  \hline
\end{tabular}
\end{center}
\clearpage

\begin{center}
$(\hF_1(\hF_1)^{Z_{22'}^{-1}Z_{66'}^{-1}})^{Z^2_{1',22'}Z^{-2}_{5',66'}}$ 
\end{center}
\begin{center}
\begin{tabular}[H]{|l|c|c|}\hline
$\vcenter{\hbox{\epsfbox{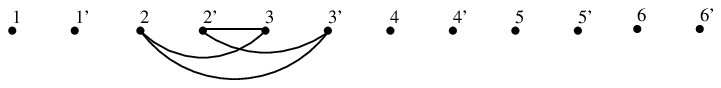}}}$  & $\rho^j_3 \rho^i_2 z_{2'3} \rho^{-i}_2 \rho^{-j}_3$ & 3\\
    \hline
$\vcenter{\hbox{\epsfbox{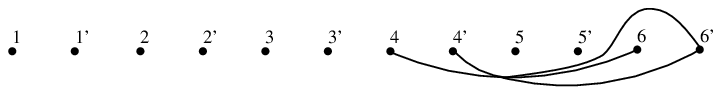}}}$  & $\rho^j_6 \rho^i_4  \underline{z}_{46} \rho^{-i}_4 \rho^{-j}_6$ & 3\\
    \hline
$\vcenter{\hbox{\epsfbox{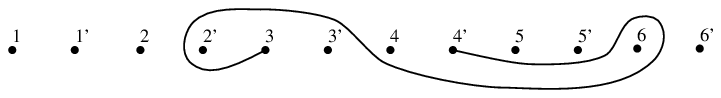}}}$  & $z_{34'}$ & 1\\
   \hline
$\vcenter{\hbox{\epsfbox{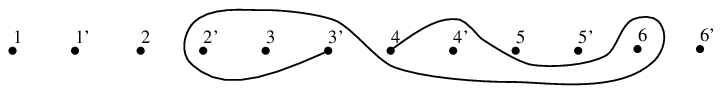}}}$   & $z_{3'4}$ & 1\\
    \hline
$\vcenter{\hbox{\epsfbox{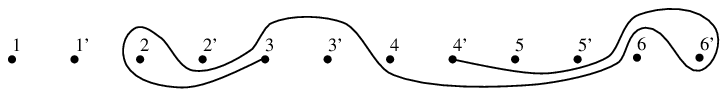}}}$  & $z_{34'}^{\rho^{-1}}$ & 1 \\  \hline
$\vcenter{\hbox{\epsfbox{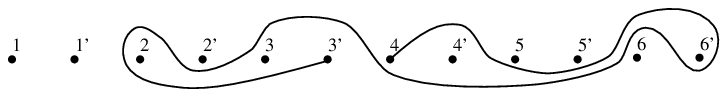}}}$  & $z_{3'4}^{\rho^{-1}}$ & 1 \\  \hline
$\vcenter{\hbox{\epsfbox{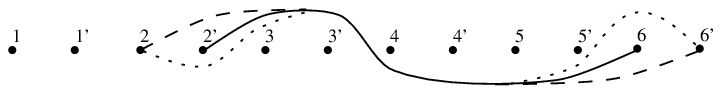}}}$   & $\rho^i_6 \rho^i_2  
\tilde{z}_{2'6} \rho^{-i}_2 \rho^{-i}_6$ & 2\\
    \hline
$\vcenter{\hbox{\epsfbox{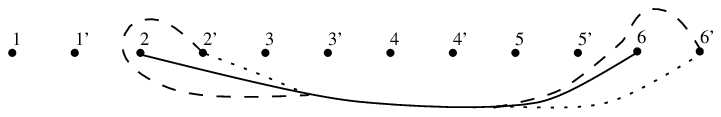}}}$   & $\rho^i_6 \rho^i_2  \underline{z}_{26} \rho^{-i}_2 \rho^{-i}_6$ & 2\\
    \hline
\end{tabular}
\end{center}

\clearpage
\noindent
$H_{V_7} = 
(\uZ^2_{2'i})_{i=33',4,4',66'} \cdot 
Z^3_{11',2} \cdot 
(\uZ^3_{2',55'})^{Z^2_{2',33'}} \cdot
(\uZ^{2}_{2i})^{Z_{11',2}}_{i = 33',4,4',66'}
 \cdot 
(Z_{22'})^{Z^2_{11',2}\bZ^2_{2',55'}Z^2_{44',55'}} \cdot
Z^3_{33',4} \cdot
(\uZ^2_{i4'})^{Z^{2}_{11',2}}_{i=11',55'} \cdot
\bZ^3_{4',66'} \cdot
(\uZ^2_{i4})^{Z^2_{33',4}Z^2_{11',2}}_{i = 11',55'} \cdot
(Z_{44'})^{\bZ^2_{4',66'}Z^2_{33',4}}
\cdot
\left(\hF_1(\hF_1)^{Z_{11'}^{-1}Z_{33'}^{-1}}\right)^{\uZ^{-2}_{4,55'}
 \uZ^{-2}_{4,66'}Z^2_{11',2}}$.\\
\begin{center}
\begin{tabular}[hb]{|l|p{3in}|l|c|}\hline
& The paths/complex conjugates & The braids & The exponent of braids \\ [-.2cm]
\hline
(1) & $\vcenter{\hbox{\epsfbox{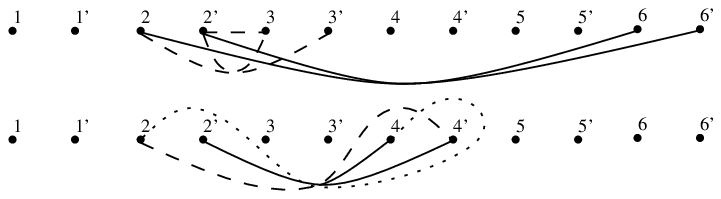}}}$  & 
$\begin{array}{lll}\rho^j_m \rho^i_2 
\underline{ z}_{2'm} \rho^{-i}_2 \rho^{-j}_m \\[-.2cm]
m= 3,6 \\ [-.2cm] 
\rho^i_4 \rho^i_2 
\underline{ z}_{2'4} \rho^{-i}_2 \rho^{-i}_4\end{array}$ & 2 \\ 
    \hline
(2) & $\vcenter{\hbox{\epsfbox{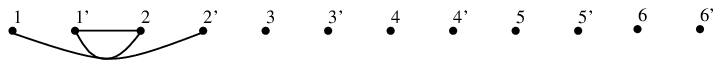}}}$ & $\rho^j_2 \rho^i_1 \underline{ z}_{12} \rho^{-i}_1 \rho^{-j}_2$ & 3  \\ 
   \hline
(3) & $\vcenter{\hbox{\epsfbox{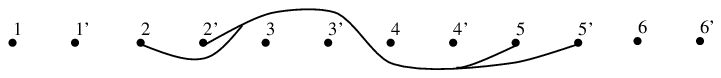}}}$ & $ \rho^j_5 \rho^i_2 \tilde{z}_{2'5} \rho^{-i}_2 \rho^{-j}_5$ & 3  
\\  
   \hline
(4) & $\vcenter{\hbox{\epsfbox{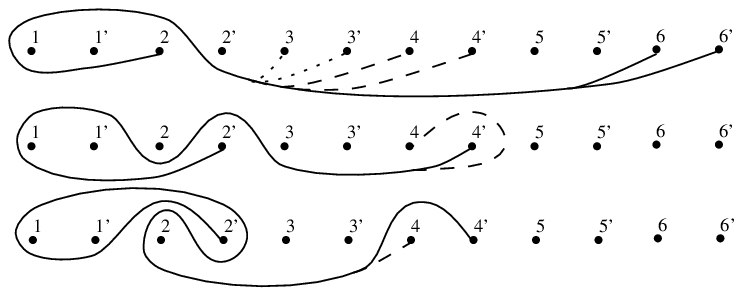}}} $ & 
$\begin{array}{lll}\rho^j_m \rho^i_2 \tilde{z}_{2m} \rho^{-i}_2 \rho^{-j}_m \\
  [-.2cm] 
   m = 3,6 \\  [-.2cm] 
\rho^i_4 \rho^i_2 \tilde{z}_{24} \rho^{-i}_2 \rho^{-i}_4\end{array}$ & 2
\\ \hline
(5) & $\vcenter{\hbox{\epsfbox{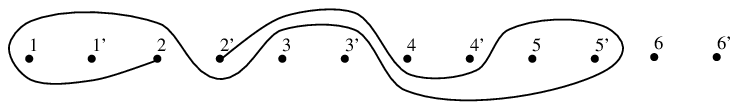}}}$  & $\tilde{z}_{22'} $ & 1\\
    \hline
(6) & $\vcenter{\hbox{\epsfbox{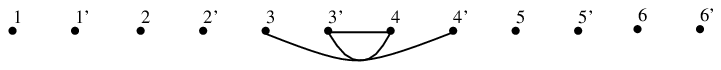}}}$  & $\rho^j_4 \rho^i_3 \underline{ z}_{34} \rho^{-i}_3 \rho^{-j}_4    $ & 3  \\ 
   \hline
(7) & $\vcenter{\hbox{\epsfbox{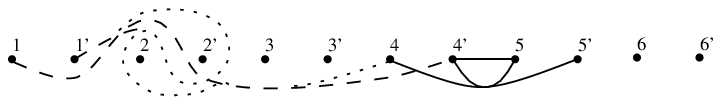}}}$  & $
\begin{array}{ll}\rho^j_m \rho^i_4 z_{4'm}\rho_4^{-i} \rho_m^{-j} \\
  [-.2cm] 
m = 1,5\end{array}$ & 2 \\  \hline
(8) & $\vcenter{\hbox{\epsfbox{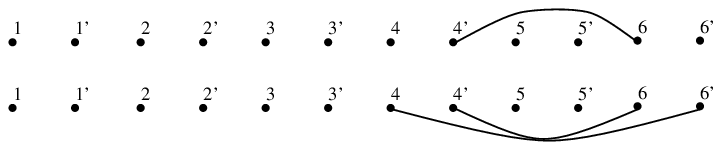}}}$  & $\rho^j_6 \rho^i_4 \underline{ z}_{4'6} \rho^{-i}_4 \rho^{-j}_6$ & 3 \\ 
    \hline
(9) & $\vcenter{\hbox{\epsfbox{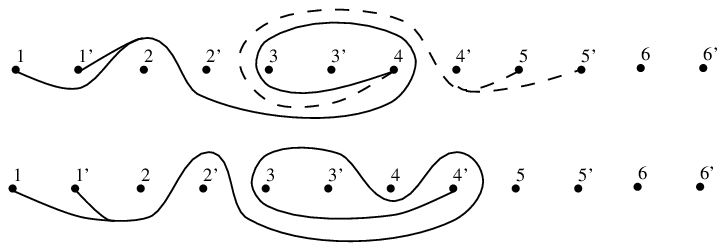}}}$  & $ \begin{array}{ll}\rho^j_4 \rho^i_m \tilde{z}_{m4} \rho^{-i}_m \rho^{-j}_4 \\
  [-.2cm] 
  m = 1,5\end{array}$ & 2 \\  \hline
(10) & $\vcenter{\hbox{\epsfbox{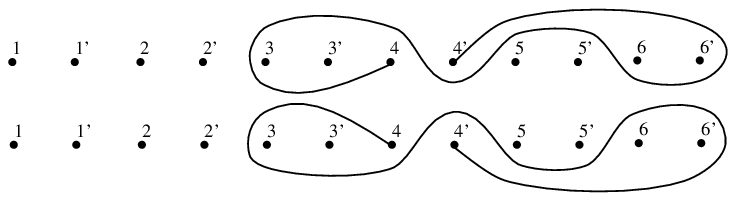}}}$ & $ \tilde{z}_{44'}$ & 1 \\ 
   \hline
\end{tabular}
\end{center}
\clearpage
\begin{center}
$(\hF_1 \cdot (\hF_1)^{Z_{11'}^{-1}Z_{33'}^{-1}})^{\uZ^{-2}_{4,55'}
 \uZ^{-2}_{4,66'}Z^2_{11',2}}$ 
\end{center}
\begin{center}
\begin{tabular}[H]{|l|c|c|}\hline
$\vcenter{\hbox{\epsfbox{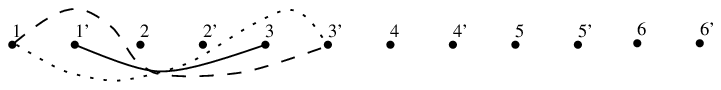}}}$ & $\rho^i_3 \rho^i_1 \underline{ z}_{1'3} \rho^{-i}_1 \rho^{-i}_3$ & 2\\
    \hline
$\vcenter{\hbox{\epsfbox{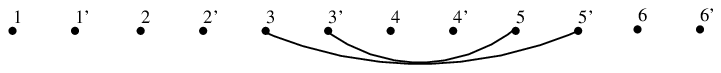}}}$  & $\rho^j_5 \rho^i_3  \underline{z}_{3'5} \rho^{-i}_3 \rho^{-j}_5$ & 3\\
    \hline
$\vcenter{\hbox{\epsfbox{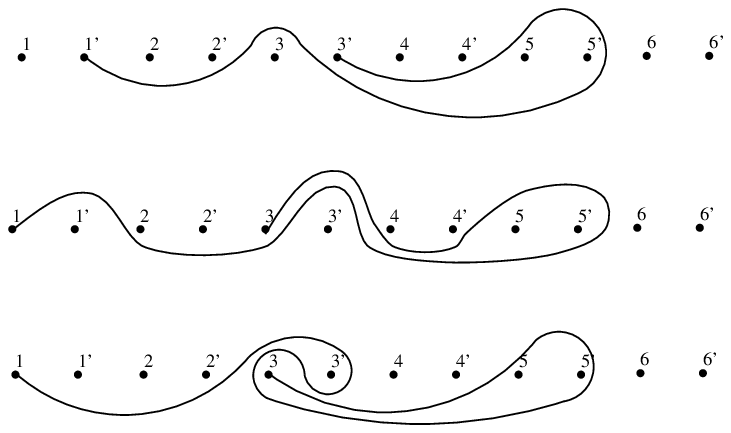}}}$ & $\rho^i_3 \rho^i_1 \tilde{z}_{1'3'} \rho^{-i}_1 \rho^{-i}_3$ & 2\\
  \hline
$\vcenter{\hbox{\epsfbox{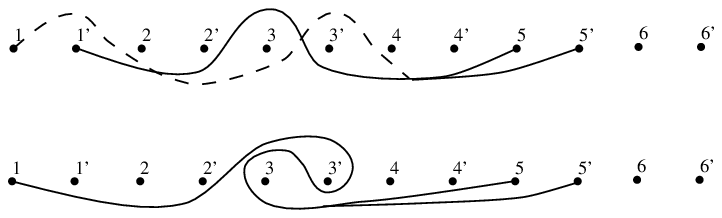}}}$ & $\rho^j_5 \rho^i_1 \tilde{z}_{1'5} \rho^{-i}_1 \rho^{-j}_5$ & 3\\
  \hline
$\vcenter{\hbox{\epsfbox{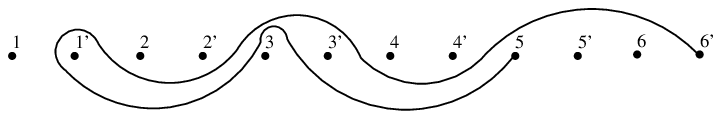}}}$ & $z_{56'}$ & 1\\
    \hline
$\vcenter{\hbox{\epsfbox{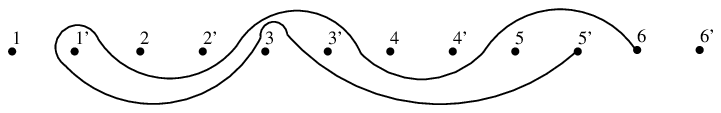}}}$  & $z_{5'6}$ & 1\\
     \hline
$\vcenter{\hbox{\epsfbox{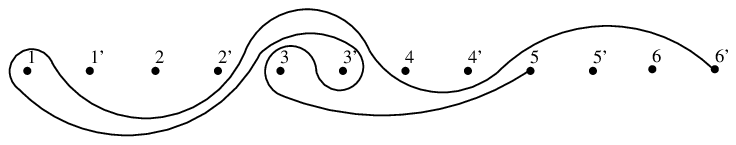}}}$  & $z_{56'}^{\rho^{-1}}$ & 1 \\
  \hline
$\vcenter{\hbox{\epsfbox{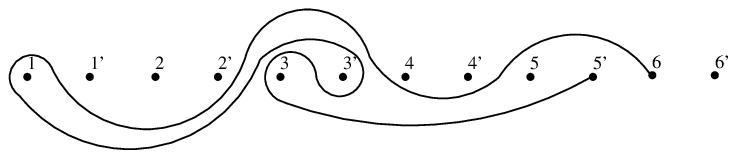}}}$  & $z_{5'6}^{\rho^{-1}}$ & 1 \\   \hline
\end{tabular}
\end{center}

\clearpage
\noindent
$H_{V_8} = Z^3_{11',2} \cdot 
(\uZ^2_{2'i})_{i=3,3',44',55'} \cdot 
(\uZ^3_{2',66'})^{Z^{-2}_{55',66'}}  \cdot 
(\uZ^2_{2'i})^{\uZ^{-2}_{i,66'}}_{i = 3,3',44',55'}
 \cdot 
(Z_{22'})^{\uZ^2_{2',66'}\uZ^2_{2',55'}Z^2_{11',2}} \cdot
Z^3_{3',44'} \cdot
(\uZ^2_{3'i})^{Z^{2}_{3',44'}Z^2_{11',2}}_{i = 11',66'} \cdot
\bZ^3_{3',55'} \cdot
(\uZ^2_{3i})^{Z^2_{11',2}}_{i = 11',66'} \cdot
(Z_{33'})^{\uZ^2_{3',55'}Z^2_{3',44'}}
\cdot
\left(\hF_1(\hF_1)^{Z_{11'}^{-1}Z_{44'}^{-1}}\right)^{\uZ^{2}_{11',3'}
 Z^{2}_{3',44'}Z^2_{11',2}}$. \\
\begin{center}
\begin{tabular}[hb]{|l|p{3in}|l|c|}\hline
& The paths/complex conjugates & The braids & The exponent of braids \\ [-.2cm]  \hline
(1) & $\vcenter{\hbox{\epsfbox{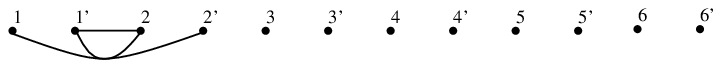}}}$ & $\rho^j_2 \rho^i_1 \underline{z}_{12} \rho^{-i}_1 \rho^{-j}_2$ & 3\\
   \hline
(2) & $\vcenter{\hbox{\epsfbox{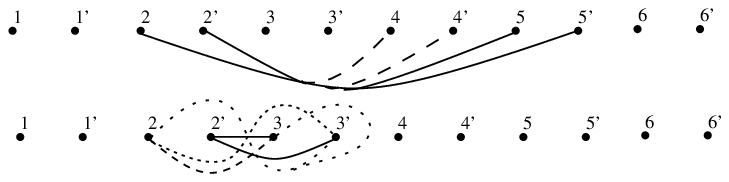}}}$ & $\begin{array}{lll}\rho^j_m \rho^i_2 \underline{z}_{2'm} \rho^{-i}_2 \rho^{-j}_m \\ [-.2cm]
 m=4,5  \\ [-.2cm]
\rho^i_3 \rho^i_2 \underline{z}_{2'3} \rho^{-i}_2 \rho^{-i}_3 \end{array}$ & 2 
\\   \hline
(3) & $\vcenter{\hbox{\epsfbox{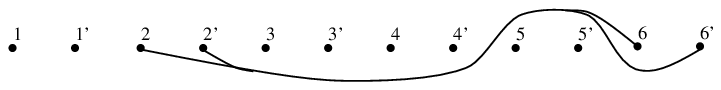}}}$ & $ \rho^j_6 \rho^i_2 \tilde{z}_{2'6} \rho^{-i}_2 \rho^{-j}_6$ & 3  
\\   \hline
(4) & $\vcenter{\hbox{\epsfbox{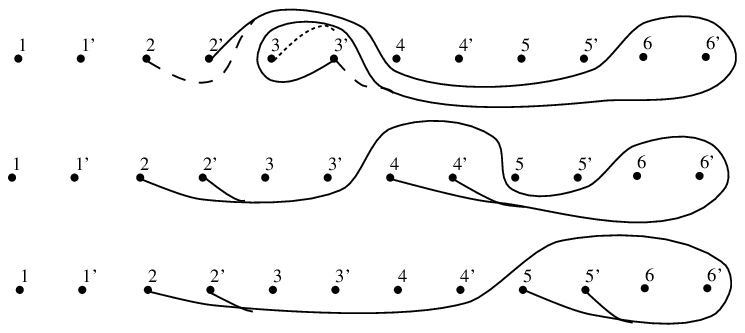}}}$ & $\begin{array}{lll}\rho^j_m \rho^i_2 \tilde{z}_{2'm} \rho^{-i}_2 \rho^{-j}_m \\ [-.2cm]
 m=4,5 \\ [.2cm]
\rho^i_3 \rho^i_2 \tilde{z}_{23} \rho^{-i}_2 \rho^{-i}_3 \end{array}$ & 2 \\
    \hline
(5) & $\vcenter{\hbox{\epsfbox{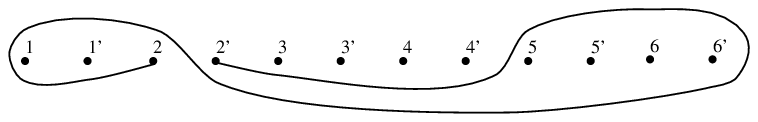}}}$ & $\tilde{z}_{22'} $ & 1 \\  
   \hline
(6) & $\vcenter{\hbox{\epsfbox{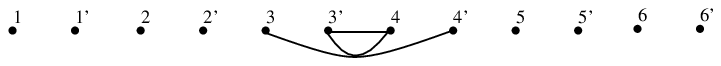}}}$ & $ \rho^j_4 \rho^i_3 \underline{z}_{3'4} \rho^{-i}_3 \rho^{-j}_4$ & 3  
\\  \hline
(7) & $\vcenter{\hbox{\epsfbox{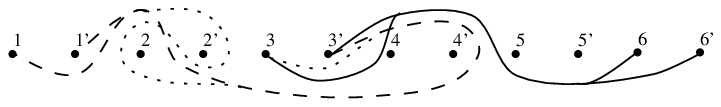}}}$  & $\rho^j_m \rho^i_3 \tilde{z}_{3'm}\rho^{-i}_3 \rho_m^{-j} $ & 2\\
[-.4cm]
 & & $m=1,6$ &\\    \hline
(8) & $\vcenter{\hbox{\epsfbox{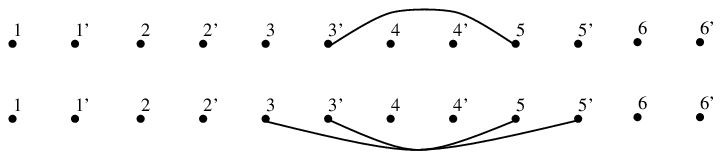}}}$ & $\rho^j_5 \rho^i_3 \underline{z}_{3'5} \rho^{-i}_3 \rho^{-j}_5$ & 3  \\ 
   \hline
(9) &  $\vcenter{\hbox{\epsfbox{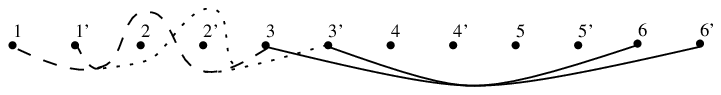}}}$ & $\rho^j_m \rho^i_3 \tilde{ z}_{3m} \rho^{-i}_3 \rho^{-j}_m$ & 2 \\ 
[-.4cm]
 & & $m=1,6$ &\\     \hline
(10) & $\vcenter{\hbox{\epsfbox{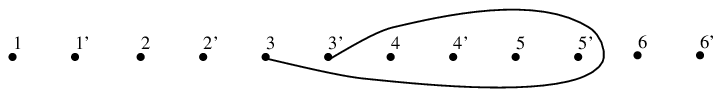}}}$ & $\tilde{z}_{33'}$ & 1 \\ 
    \hline
\end{tabular}
\end{center}
\clearpage
\begin{center}
$(\hF_1 \cdot (\hF_1)^{Z_{11'}^{-1}Z_{44'}^{-1}})^{\uZ^{2}_{11',3'}
 Z^{2}_{3',44'}Z^2_{11',2}}$ 
\end{center}
\begin{center}
\begin{tabular}[H]{|l|c|c|}\hline
$\vcenter{\hbox{\epsfbox{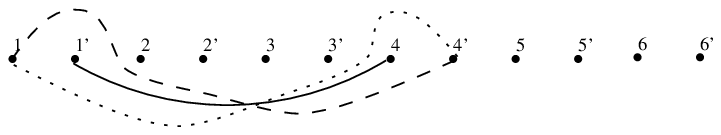}}}$  & $\rho^i_4 \rho^i_1 \underline{ z}_{1'4} \rho^{-i}_1 \rho^{-i}_4$ & 2\\
   \hline
$\vcenter{\hbox{\epsfbox{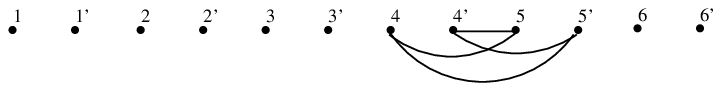}}}$   & $\rho^j_5 \rho^i_4  \underline{z}_{4'5} \rho^{-i}_4 \rho^{-j}_5$ & 3\\
   \hline
$\vcenter{\hbox{\epsfbox{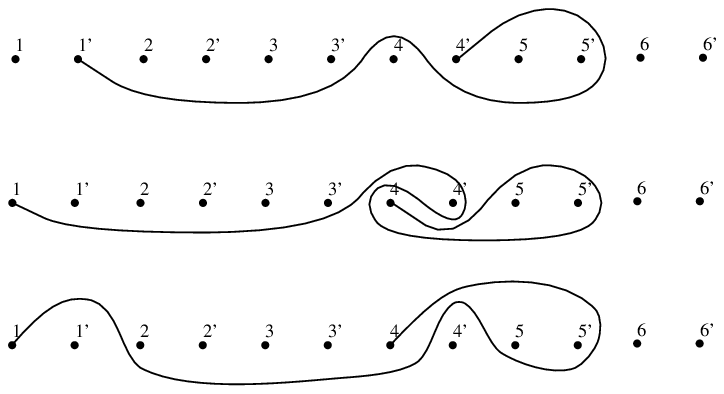}}}$  & $\rho^i_4 \rho^i_1 \tilde{z}_{1'4'} \rho^{-i}_1 \rho^{-i}_4$ & 2\\
  \hline
$\vcenter{\hbox{\epsfbox{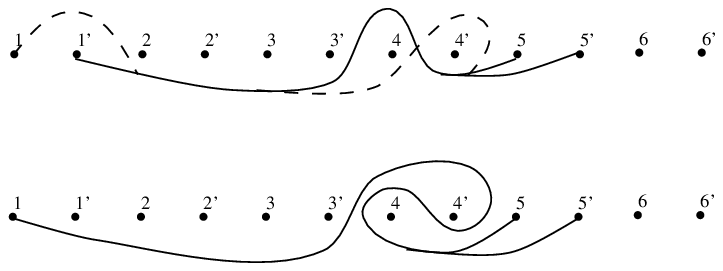}}}$  & $\rho^j_5 \rho^i_1 \tilde{z}_{1'5} \rho^{-i}_1 \rho^{-j}_5$ & 3\\
  \hline
$\vcenter{\hbox{\epsfbox{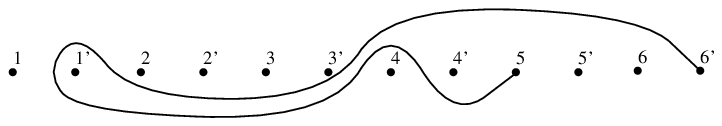}}}$  & $z_{56'}$ & 1\\
    \hline
$\vcenter{\hbox{\epsfbox{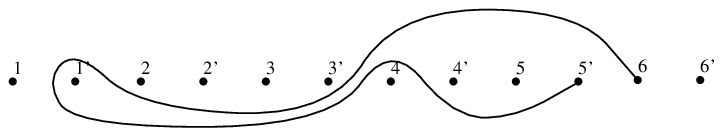}}}$   & $z_{5'6}$ & 1\\
    \hline
$\vcenter{\hbox{\epsfbox{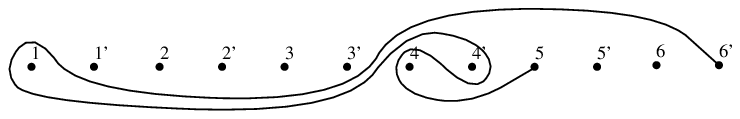}}}$  & $z_{56'}^{\rho^{-1}}$ & 1 \\     \hline
$\vcenter{\hbox{\epsfbox{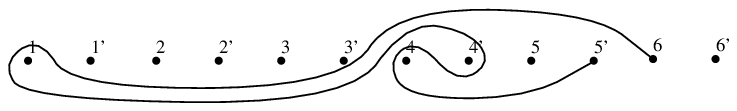}}}$  & $z_{5'6}^{\rho^{-1}}$ & 1 \\  \hline
\end{tabular}
\end{center}
\medskip
\noindent
$H_{V_9} = Z^3_{1',22'} \cdot 
(\uZ^2_{1'i})^{Z^2_{1',22'}}_{i=3,3',44',66'}
 \cdot 
(\uZ^3_{1',55'})^{Z^2_{1',22'}}  \cdot
(\uZ^2_{1'i})^{\uZ^2_{1',55'}Z^2_{1',22'}}_{i = 3,3',44',66'} 
\cdot
(Z_{11'})^{\uZ^2_{1',66'}\uZ^2_{1',55'}Z^2_{1',22'}} \cdot
Z^3_{3',44'} \cdot
(\uZ^2_{3i})^{Z^2_{1',22'}}_{i=22',55'}
\cdot
(\uZ^2_{i3'})^{Z^2_{3',44'}\uZ^2_{i3}Z^2_{1',22'}}_{i = 22',55'} \cdot
(\uZ^3_{3',66'})^{Z^2_{3',44'}} \cdot
(Z_{33'})^{\uZ^2_{3',66'}Z^2_{3',44'}} \cdot
(\hF_1(\hF_1)^{Z_{22'}^{-1}Z_{44'}^{-1}})^{Z^2_{22',33'}Z^{2}_{3',44'} Z^{2}_{1',22'}}$. \\
\begin{center}
\begin{tabular}[hb]{|l|p{3in}|l|c|}\hline
& The paths/complex conjugates & The braids & The exponent of braids \\ [-.2cm]
\hline
(1) & $\vcenter{\hbox{\epsfbox{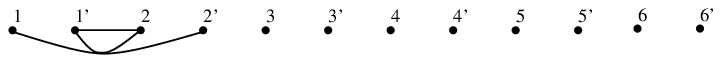}}}$ & $\rho^j_2 \rho^i_1  z_{1'2} \rho^{-i}_1 \rho^{-j}_2$ & 3\\
    \hline
(2) & $\vcenter{\hbox{\epsfbox{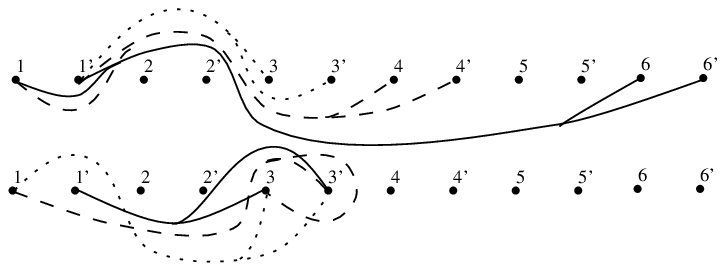}}}$ & $\begin{array}{lll}\rho^j_m \rho^i_1 \tilde{z}_{1'm} \rho^{-i}_1 \rho^{-j}_m \\ [-.2cm]
 m=4,6  \\ [-.2cm]
\rho^i_3 \rho^i_1 \underline{z}_{1'3} \rho^{-i}_1 \rho^{-i}_3 \end{array}$ & 2 
\\   \hline
(3) & $\vcenter{\hbox{\epsfbox{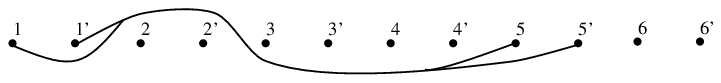}}}$ & $ \rho^j_5 \rho^i_1 \tilde{ z}_{
1'5} \rho^{-i}_1 \rho^{-j}_5$ & 3 
\\  
   \hline
(4) & $\vcenter{\hbox{\epsfbox{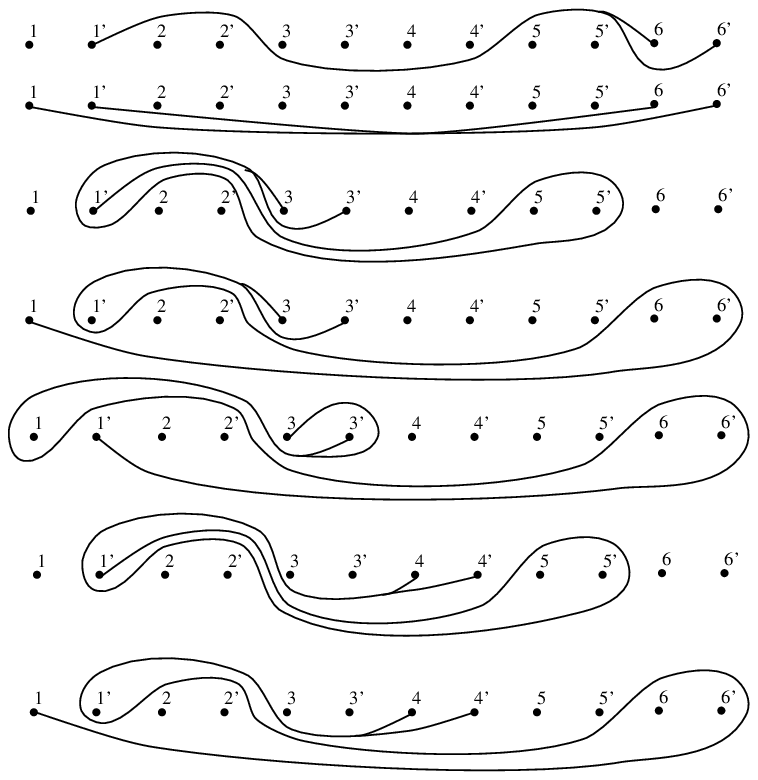}}}$ & $\begin{array}{lll}\rho^j_m \rho^i_1 \tilde{z}_{1'm} \rho^{-i}_1 \rho^{-j}_m \\ [-.2cm]
 m=4,6  \\ [-.2cm]
\rho^i_3 \rho^i_1 \tilde{z}_{13} \rho^{-i}_1 \rho^{-i}_3 \end{array}$ & 2 
\\   \hline
(5) & $\vcenter{\hbox{\epsfbox{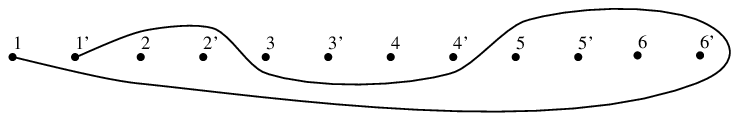}}}$ & $ \tilde{z}_{11'}$ & 1\\
   \hline
(6) & $\vcenter{\hbox{\epsfbox{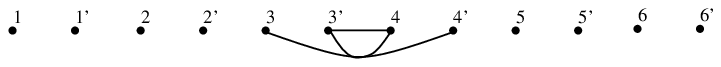}}}$  & $  \rho^j_4 \rho^i_3 z_{3'4} \rho^{-i}_3 \rho^{-j}_4 $ & 3  \\ 
   \hline
(7) & $\vcenter{\hbox{\epsfbox{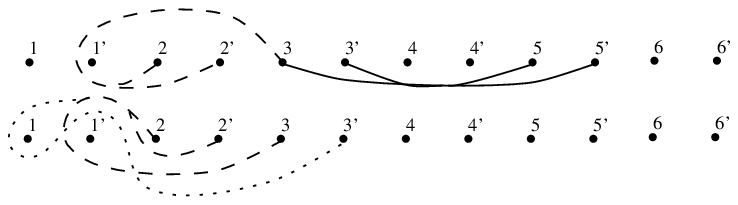}}}$ & $ \begin{array}{ll}\rho^j_m \rho^i_3  \tilde{z}_{3m} \rho^{-i}_3 \rho^{-j}_m \\ [-.2cm]
m = 2,5 \end{array}$ & 2 \\  
   \hline
(8) & $\vcenter{\hbox{\epsfbox{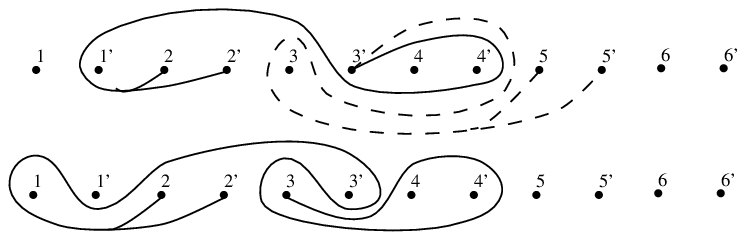}}}$ &  $\begin{array}{ll}\rho^j_m \rho^i_3  \tilde{z}_{3'm} \rho^{-i}_3 \rho^{-j}_m \\ [-.2cm]
m = 2,5 \end{array}$ & 2 \\  
   \hline
(9) & $\vcenter{\hbox{\epsfbox{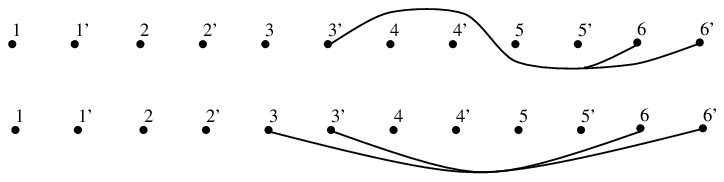}}}$ & $ \rho^j_6 \rho^i_3 \tilde{z}_{3'6} \rho^{-i}_3 \rho^{-j}_6$ & 3 \\ 
    \hline
(10) & $\vcenter{\hbox{\epsfbox{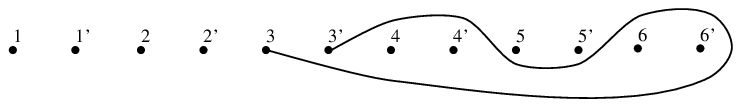}}}$ & $ \tilde{z}_{33'}$ & 1 \\ 
   \hline
\end{tabular}
\end{center}
\clearpage
\begin{center}
$(\hF_1 \cdot (\hF_1)^{Z_{22'}^{-1}Z_{44'}^{-1}})^{Z^2_{22',33'}Z^{2}_{3',44'} Z^{2}_{1',22'}}$ 
\end{center}


\begin{center}
\begin{tabular}[H]{|l|c|c|}\hline
$\vcenter{\hbox{\epsfbox{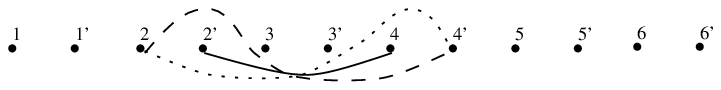}}}$ & $\rho^i_4 \rho^i_2 \underline{ z}_{2'4} \rho^{-i}_2 \rho^{-i}_4$ & 2\\
    \hline
$\vcenter{\hbox{\epsfbox{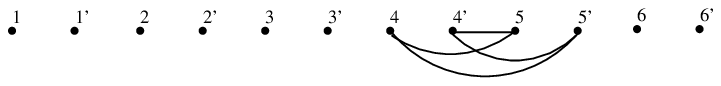}}}$  & $\rho^j_5 \rho^i_4  z_{4'5} \rho^{-i}_4 \rho^{-j}_5$ & 3\\
    \hline
$\vcenter{\hbox{\epsfbox{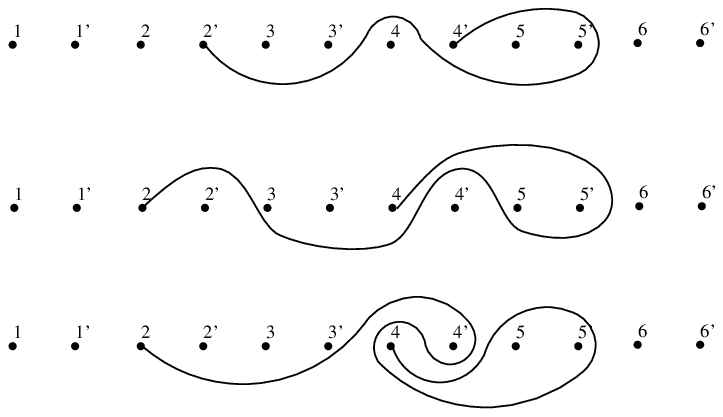}}}$ & $\rho^i_4 \rho^i_2 \tilde{z}_{2'4'} \rho^{-i}_2 \rho^{-i}_4$ & 2\\
  \hline
$\vcenter{\hbox{\epsfbox{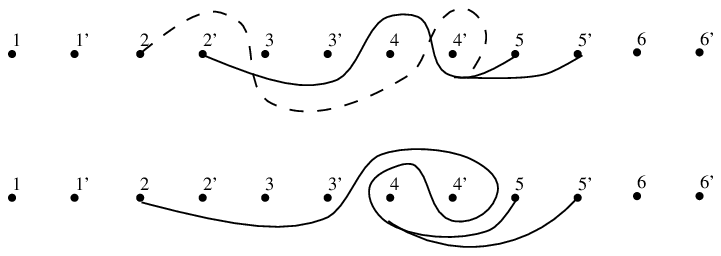}}}$ & $\rho^j_5 \rho^i_2 \tilde{z}_{2'5}
\rho^{-i}_2 \rho^{-j}_5 $ & 3\\
  \hline
$\vcenter{\hbox{\epsfbox{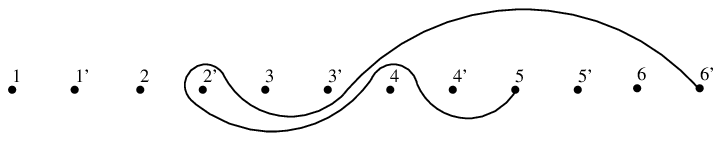}}}$ & $z_{56'}$ & 1\\
     \hline
$\vcenter{\hbox{\epsfbox{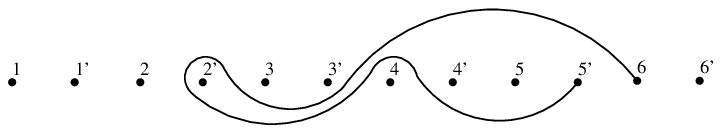}}}$  & $z_{5'6}$ & 1\\
     \hline
$\vcenter{\hbox{\epsfbox{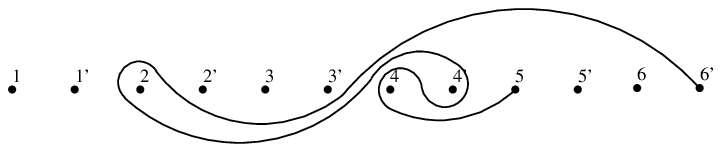}}}$  & $z_{56'}^{\rho^{-1}}$ & 1 \\   \hline
$\vcenter{\hbox{\epsfbox{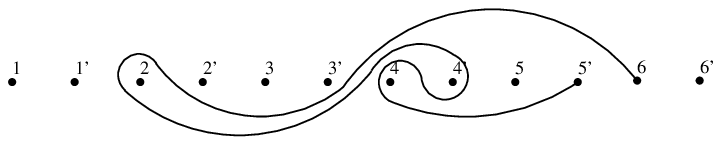}}}$  & $z_{5'6}^{\rho^{-1}}$ & 1 \\
   \hline
\end{tabular}
\end{center}

\subsection{Properties of $\Delta^2_{54}$}
\indent

Now we are checking if the braids we have found are the only ones.For that, we assume that $\Delta^2_{54} = \prod\limits^9_{i=1} C_i H_{V_i}\prod b_i$ be the braid monodromy factorization. $b_i$ are factors corresponding to singularities that are not covered by $\prod\limits^9_{i=1} C_i H_{V_i}$, and each
  $b_i$ is of the form $Y^{t_i}_i$, $Y_i$
  is a positive halftwist and $0 \leq t_i \leq 3$.  We compute degrees of all factors involved.
By [6, Theorem 8], 
deg $\Delta^2_{54} = 54 \cdot 53 = 2862$.

On the other hand, $\prod\limits^9_{i=1} C_i$ compounds 864 factors,
so deg$\prod\limits^9_{i=1} C_i = 864 \cdot 2 = 1728$. In the computations of each
$H_{V_i} , 1 \leq i \leq 9$, deg$\hat{F}_1 = \mbox{deg}(\hat{F}_1)^{\rho^{-1}} = 24$, so
 deg$\hat{F}_1(\hat{F}_1)^{\rho^{-1}} = 48$.  The factors outside
$\hat{F}_1(\hat{F}_1)^{\rho^{-1}}$ in each $H_{V_i}$ are: 20
degree two factors,
 12 degee three factors, 2 degree one factor.
 The sum of the factors' degrees resulting from each $H_{V_i}$ is $20
 \cdot 2 + 12 \cdot 3 + 2 \cdot 1 = 78$. Thus deg$H_{V_i} = 78 + 48 = 126$ for $1
  \leq i \leq 9$, and therefore deg$\prod\limits^9_{i=1}H_{V_i} = 126 \cdot 9 = 1134$.
   Finally deg$\prod\limits^9_{i=1}C_i H_{V_i} = 1134 + 1728 = 2862$.

Therefore, deg$\prod\limits^9_{i=1}b_i = 2862 \cdot 2862^{-1} = 1$.
 Since $\forall i, b_i$ is a positive exponent of a positive halftwist,
  we get $b_i = 1 \ \forall i$.

Finally, $\Delta^2_{54} = \prod\limits^9_{i=1} C_i H_{V_i}$.

In the following lemma we prove the Invariance Property for
$C_i$.

\begin{lemma} {\bf Complex Invariance of $C_i$.}\label{invCi}
For every $i , i = 1, \cdots , 9$ and $\forall m_j \in \Z \ , \ 1 \leq j \leq 27 \ , C_i$ is invariant under $\prod\limits^{27}_{j=1} Z^{m_j}_{jj'}$.
\end{lemma}

\begin{proof}
We apply Invariance Rule II and 
Invariance Remark (iv) [6, Subsection 7.3] 
on each of the factors of $C_i$ of the form $Y^{2}_{ii',jj'}$.
\end{proof}

We can summarize the properties of $\Delta^2_{54}$ in the following theorems.

\begin{theorem} {\bf Invariance Theorem for $\Delta^2_{54}$.}\label{invdel}
Let $\rho = \prod\limits^{27}_{j=1} \rho_{m_j}$ for $\rho_{m_j} = Z^{m_j}_{jj'}$ and $1 \leq j \leq 27$.  Then $\Delta^2_{54}$ is invariant under $\rho$ for every $m_j \in \Z$.
\end{theorem}

\begin{proof}  
By Invariance Properties of
$H_{V_i}$ ([4 , Sections 3.1 - 3.9]) and by Lemma \ref{invCi}.
\end{proof}

\begin{theorem} {\bf Complex Conjugation Theorem.}\label{complex}
$\Delta^2_{54}$ is invariant under complex conjugation.
\end{theorem}

\begin{proof} 
The same proof from Lemma 19, [22].
\end{proof}

\subsection{Consequences from the Invariance Theorems}
\indent

Recall that $S^{(0)} = S$ is the regenerated branch curve.  Every
factor of a braid monodromy factorization $\Delta^2_{54}$ induces
a relation on $\pi_1(\C^2 - S,M)$. The Invariance Properties are
an essential addition to the \vK \ Theorem, since we get more
relations in   $\pi_1(\C^2 - S,M)$.

\begin{theorem}\label{th:161} {\bf  [6, Sction 3.10]}.   
Let $S, \varphi, B_{54}$ be as above. If a sub-factorization 
$\prodl^r_{i=s} Z_i$ of $\Delta^2_{54}$ is invariant 
under any element $h$,  and $\prodl^r_{i=s} Z_i$ induces a relation 
$\G_{i_1} \cdot ... \cdot \G_{i_t}$ on  $\pi_1(\C^2 - S,M)$ 
via the \vK \ method, then 
$(\G_{i_1})_h  \cdot  ... \cdot (\G_{i_t})_h$ is also a relation.
\end{theorem}

\begin{corollary}\label{infinite} 
If $R$ is any relation in
 $\pi_1(\C^2 - S,M)$, then $R_\rho$ is also a relation in  $\pi_1(\C^2 - S,M)$, 
where $R_\rho$ is
the relation induced from  $R$ by replacing $\G_j$ and $\G_{j'}$ with
$(\G_j)_{\rho^{m_j}_j}$ and $(\G_{j'})_{\rho^{m_j}_j}$ 
respectively for $1 \leq j \leq 27$,
where $\rho = \prodl^{27}_{j=1} \rho_{m_j}$ 
and $\rho_{m_j} = Z^{m_j}_{jj'}$.
\end{corollary}

\begin{proof} 
By Theorem \ref{invdel}.
\end{proof}

\begin{corollary}  
$\pi_1(\C^2 - S, M)$
satisfies all the relations induced in $R(\Delta^2_{54})$, all the
relations induced from $(R(\Delta^2_{54}))_{\rho}$  and all
relations induced from the complex conjugation of $\Delta^2_{54}$.
\end{corollary}

\begin{proof} 
By Theorems \ref{invdel} and \ref{complex}.
\end{proof}

\section{$\pi_1(\C^2 - S, M)$ and $\tilde{\pi}_1$}

\subsection{The \vK\ Theorem}\label{van}
\indent

The \vK\ Theorem induces a finite presentation of the fundamental group of complements of curves by meaning of generators and relations.

\def\suchthat{{\,:\,\,}}

We obtained  the regenerated braid monodromy factorization 
$\Delta^2_{54} = \prodl^9_{i=1} C_i H_{V_i}$.  
We have to apply the \vK\ Theorem on
the paths, which correspond to the factors in $\Delta^2_{54}$. We take any
path from $k$ to $\ell$, cut it in $M$, then towards
$k$ along the path, around $k$ and coming back the same way.
Consider $A$ as an element of the fundamental group.  Do the
same to $\ell$ to obtain $B$.  So $A, B$ are conjugations
of $\G_k$ and $\G_\ell$ respectively, see Figure 5.

\begin{figure}[h]\label{AB}
\begin{center}
\epsfbox {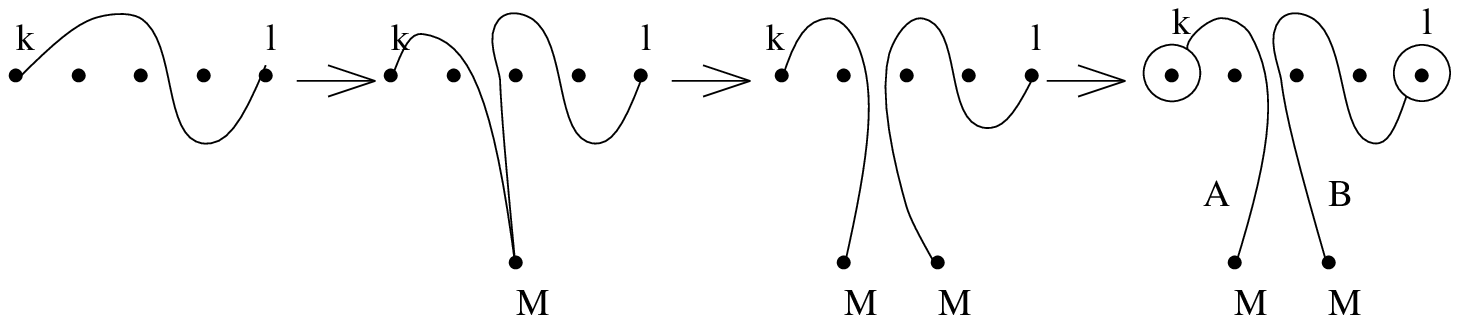}
\end{center}
\caption{}\end{figure}

\begin{lemma}\label{lm43}{\bf [40]}  $\{A, B\}$ can be extended to a 
$g$-base of $\pi_1(\C_u - S, M)$ denoted by $\{\G_j, \G_{j'}\}^{27}_{j=1}$.\\
Moreover,  there is an epimorphism $\pi_1(\C_u - S, M) \rightarrow 
\pi_1(\C^2 - S, M)$.
\end{lemma}

\begin{theorem}{\bf van Kampen for cuspidal curves.}\label{vKcuspidal}
Let $S$ be the regenerated cuspidal branch  curve, $u, M, \varphi, A,
B$ defined as above.  Let $\{\delta_i\}$ be a
$g$-base of $\pi_1(E - N, u)$.  Let $\varphi(\delta_i) =
V^{\nu_i}_i, V_i$ be a halftwist, $\nu_i = 1,2,3$.
\def\C{{\Bbb C}}
\def\P{{\Bbb P}}
Let $\set{\G_j,\G_{j'}}^{27}_{i=1}$ be a
$g$-base for $\pi_1 (\C_u - S, M)$. 
Then: $\pi_1(\C^2 - S,M)$ is generated by the images of 
$\G_j,\G_{j'}$ in $\pi_1(\C^2 - S,M)$ (denoted also by $\G_j,\G_{j'}$) 
and we get a complete set of relations from those induced
from
$\varphi (\delta_i) = V^{\nu_i}_i$, as follows (when $A, B$
are expressed in terms of $\{\G_j,\G_{j'} \})$:

(a). $A = B$, when $\nu_i = 1$;

(b). $[A,B] = 1$, when $\nu_i = 2$;

(c). $<A,B> =A B A B A B= 1$, when $\nu_i = 3$.
\end{theorem}

Recall that 
$$\pitil = {\pi_1(\C^2-S,M) \over\sg{\Gamma_j^2, \G_{j'}}}.$$

Let us apply the theorem on three examples, taken from the table corresponding to $H_{V_2}$.
We follow Figures 6, 7, 8. 
\begin{figure}[h]\label{vk1}
\epsfxsize=9cm 
\epsfysize=5.7cm 
\begin{center}
\epsfbox {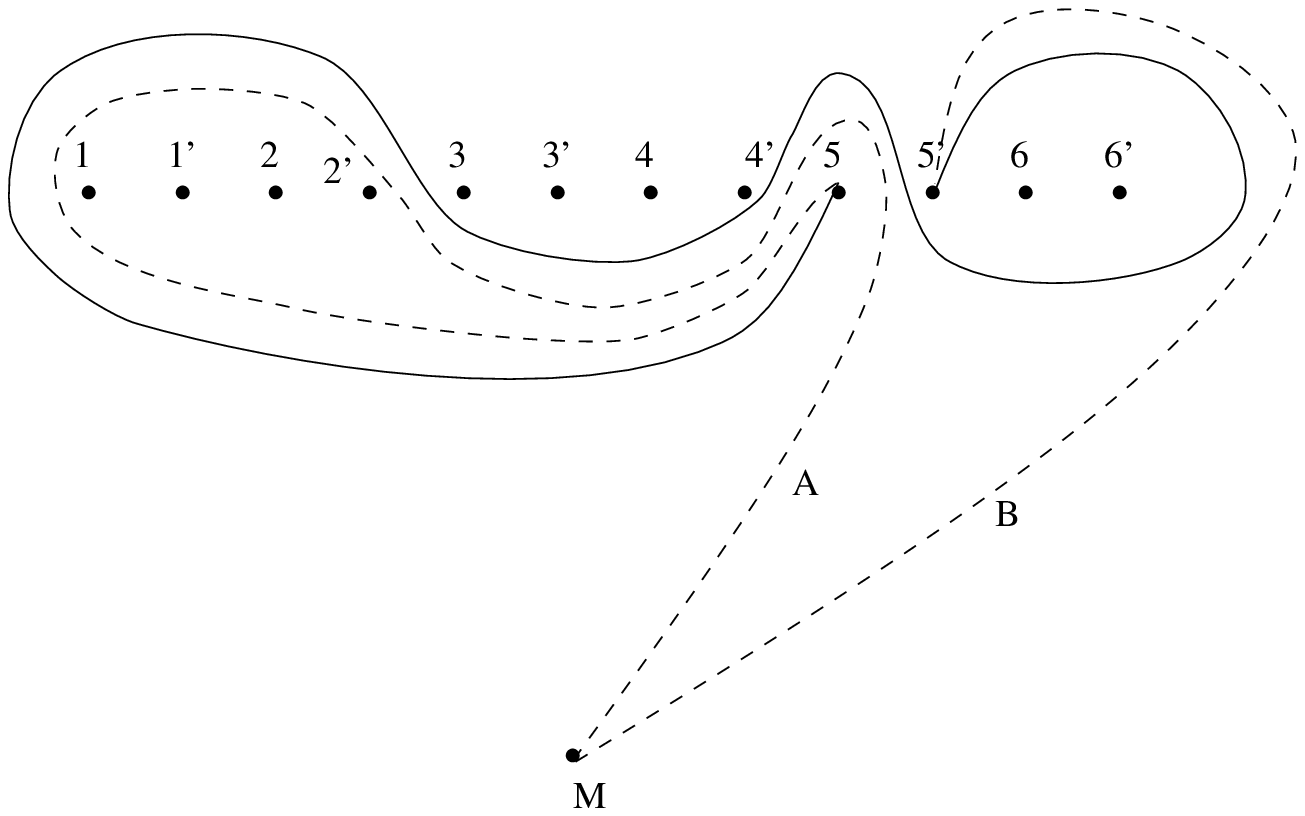}
\end{center}
\caption{}\end{figure}
\begin{figure}[h]\label{vk2}
\epsfxsize=9cm 
\epsfysize=5.7cm 
\begin{center}
\epsfbox {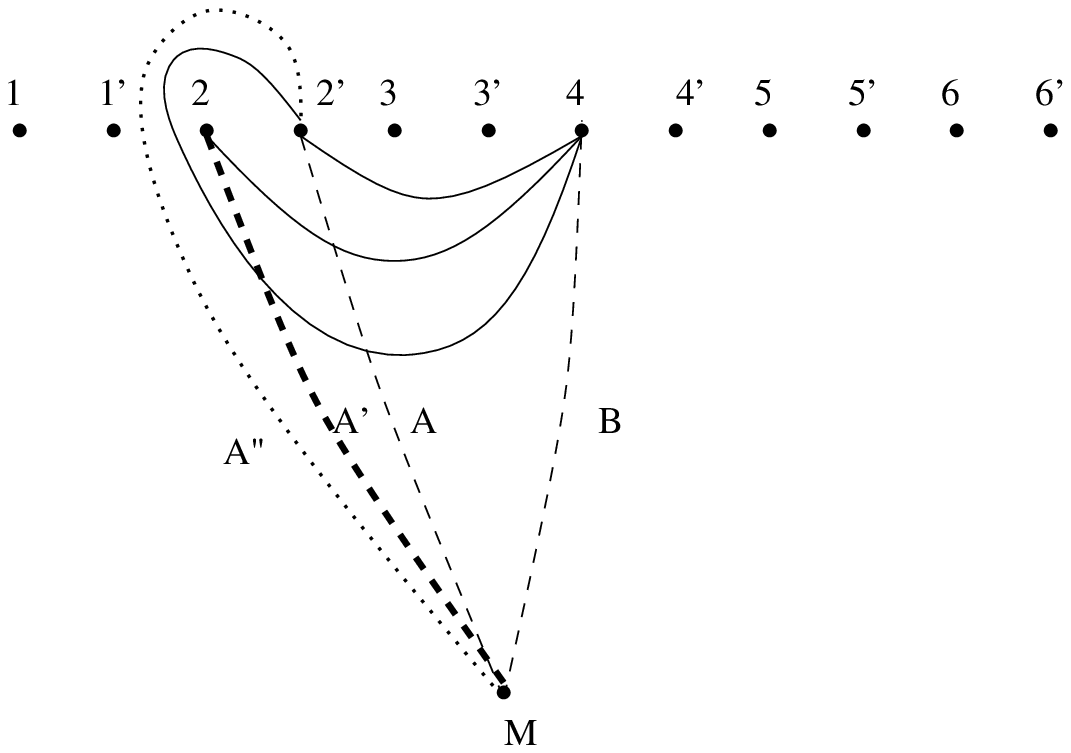}
\end{center}
\caption{}\end{figure}
\begin{figure}[h]\label{vk3}
\epsfxsize=9cm 
\epsfysize=5.7cm 
\begin{center}
\epsfbox {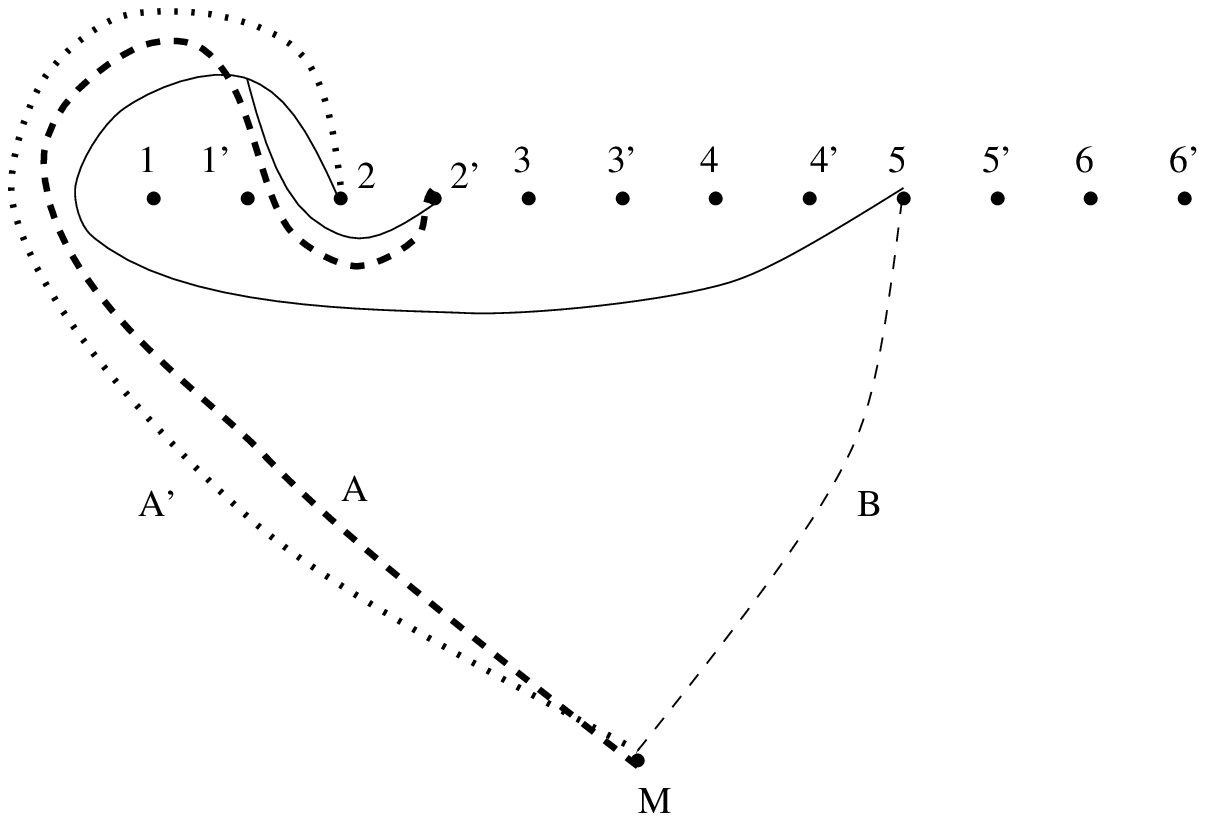}
\end{center}
\caption{}\end{figure}

Following Figure 6, $A$ and $B$ are conjugations of $\G_5$ and
$\G_{5'}$ respectively.  All factors in the conjugations are with
a positive exponent since $\G_j = \G^{-1}_j$ and $\G_{j'} =
\G^{-1}_{j'}$ in $\pitil$.

We construct $B$ by proceeding from $M$ towards $5'$ above $6'$ and 6 encircling $5'$ counterclockwise and proceeding back above 6 and $6'$.  Therefore $B = \G_{6'} \G_6 \G_{5'} \G_6 \G_{6'} = (\G_{5'})^{\G_6 \G_{6'}}$.  In a similar way, $A = \G_5 \G_{2'} \G_2 \G_{1'} \G_1 \G_5 \G_1 \G_{1'} \G_2 \G_{2'} \G_5 = (\G_5)^{\G_1 \G_{1'} \G_2 \G_{2'} \G_5}$.

This path is related to a braid which is induced from a branch point.
Therefore by Theorem \ref{vKcuspidal}, $A=B$.  The derived relation is: $(\G_5)^{\G_1 \G_{1'}} 
\G_2 \G_{2'} \G_5 = (\G_{5'})^{\G_6 \G_{6'}}$.

Following the same technique in Figure 7 we get three relations
according to three paths.  These paths represent the three braids
derived from three cusps. The relations are: $\la\G_{2'} , \G_4
\ra = 1 , \la\G_2, \G_4\ra = 1 , \la(\G_{2'})^{\G_2} , \G_4\ra =1$

Following the same technique in Figure 8, we get two relations corresponding to the two paths.  These paths represent  the two braids derived from two nodes.
The relations are: $[\G_5, (\G_2)^{\G_{1'}\G_1}] = 1$ and $[\G_5, (\G_{2'})^{\G_{1'}\G_1}] = 1$.

Since we followed the process on the point $V_2$, we continue and present the list of relations for V$_2$. We get an infinite number of relations by 
Corollary \ref{infinite}.
We present partially the infinite list for $V_2$, following the list of paths and braids involved in $H_{V_2}$ (Subsection 2.2). 
Recall that the generators involved in this certain case are $\set{\G_j, \G_{j'}}_{i=1,3,7,9,14,17}$ (see [4 , Figure 6]) which are numerated locally as 
  $\set{\G_j, \G_{j'}}_{i=1,2,3,4,5,6}$. Moreover, we add the relations 
$\G_j^2=\G_{j'}^2=1$.

\subsection{Relations $H_{V_2}$}

\indent

The following relations are induced
from the factors in $H_{V_2}$ and their paths in the first table.\\
$Relation[1]: \ \ <\G_3,\G_4>=1 \\
Relation[2]: \ \ <\G_3,\G_{4'}>=1 \\
Relation[3]: \ \ <\G_{3'},\G_4>=1 \\
Relation[4]: \ \ <\G_{3'},\G_{4'}>=1 \\
Relation[5]: \ \ [\G_1,\G_4^{(\G_{2'} \cdot \G_2)}]=1 \\
Relation[6]: \ \ [\G_1,\G_{4'}^{(\G_{2'} \cdot 
\G_2)}]=1 \\
Relation[7]: \ \ [\G_{1'},\G_4^{(\G_{2'} \cdot 
\G_2)}]=1 \\
Relation[8]: \ \ [\G_{1'},\G_{4'}^{(\G_{2'} \cdot 
\G_2)}]=1 \\
Relation[9]: \ \ [\G_4^{(\G_{2'} \cdot \G_2)},\G_5]=1 \\
Relation[10]: \ \ [\G_4^{(\G_{2'} \cdot \G_2)},\G_{5'}]=1 \\
Relation[11]: \ \ [\G_{4'}^{(\G_{2'} \cdot \G_2)},\G_{5'}]=1 \\
Relation[12]: \ \ [\G_{4'}^{(\G_{2'} \cdot \G_2)},\G_5^{\G_{5'}}]
=1 \\
Relation[13]: \ \ [\G_{4'}^{(\G_4 \cdot 
\G_{2'} \cdot \G_2)},\G_5]=1 \\
Relation[14]: \ \ [\G_{4'}^{(\G_4 \cdot \G_{2'} \cdot \G_2)},
\G_{5'}^{\G_5}]=1 \\
Relation[15]: \ \ [\G_4^{(\G_{2'} \cdot \G_2)},\G_6]=1 \\
Relation[16]: \ \ [\G_4^{(\G_{2'} \cdot \G_2)},\G_{6'}]=1 \\
Relation[17]: \ \ [\G_{4'}^{(\G_{2'} \cdot \G_2)},\G_6]=1 \\
Relation[18]: \ \ [\G_{4'}^{(\G_{2'} \cdot \G_2)},\G_{6'}]=1 \\
Relation[19]: \ \ <\G_2,\G_4>=1 \\
Relation[20]: \ \ <\G_2,\G_{4'}>=1 \\
Relation[21]: \ \ <\G_{2'},\G_4>=1 \\
Relation[22]: \ \ <\G_{2'},\G_{4'}>=1 \\
Relation[23]: \ \ [\G_1,\G_4]=1 \\
Relation[24]: \ \ [\G_1,\G_{4'}]=1 \\
Relation[25]: \ \ [\G_{1'},\G_4]=1 \\
Relation[26]: \ \ [\G_{1'},\G_{4'}]=1 \\
Relation[27]: \ \ [\G_4^{(\G_{3'} \cdot \G_3 \cdot \G_{2'} \cdot 
\G_2)}, \G_5^{\G_{5'}}]=1 \\
Relation[28]: \ \ [\G_4^{(\G_{3'} \cdot \G_3 \cdot \G_{2'} \cdot 
\G_2)},\G_{5'}]=1 \\
Relation[29]: \ \ [\G_{4'}^{(\G_4 \cdot \G_{3'} \cdot \G_3 \cdot \G_{2'} 
\cdot \G_2)},\G_5]=1 \\
Relation[30]: \ \ [\G_{4'}^{(\G_4 \cdot \G_{3'} \cdot \G_3 \cdot 
\G_{2'} \cdot \G_2)},\G_{5'}]=1 \\
Relation[31]: \ \ [\G_4^{(\G_{4'} \cdot \G_4 \cdot \G_{3'} \cdot 
\G_3 \cdot \G_{2'} \cdot \G_2)},\G_5]=1 \\
Relation[32]: \ \ [\G_4^{(\G_{4'} \cdot \G_4 \cdot \G_{3'} \cdot 
\G_3 \cdot \G_{2'} \cdot \G_2)},\G_{5'}^{\G_5}]=1 \\
Relation[33]: \ \ [\G_4^{(\G_{3'} \cdot \G_3 \cdot \G_{2'} \cdot \G_2)},\G_6]=1 \\
Relation[34]: \ \ [\G_4^{(\G_{3'} \cdot \G_3 \cdot \G_{2'} \cdot 
\G_2)},\G_{6'}]=1 \\
Relation[35]: \ \ [\G_{4'}^{(\G_{3'} \cdot \G_3 \cdot \G_{2'} \cdot \G_2)},\G_6]=1 \\
Relation[36]: \ \ [\G_{4'}^{(\G_{3'} \cdot \G_3 \cdot \G_{2'} \cdot \G_2)},\G_{6'}]=1 \\
Relation[37]: \ \ \G_4^{(\G_{3'} \cdot \G_3 \cdot \G_{2'} \cdot \G_2)}
=\G_{4'} \\
Relation[38]: \ \ <\G_5,\G_6>=1 \\
Relation[39]: \ \ <\G_5,\G_{6'}>=1 \\
Relation[40]: \ \ <\G_{5'},\G_6>=1 \\
Relation[41]: \ \ <\G_{5'},\G_{6'}>=1 \\
Relation[42]: \ \ [\G_3^{\G_4},\G_5]=1 \\
Relation[43]: \ \ [\G_3^{\G_{4'}},\G_{5'}]=1 \\
Relation[44]: \ \ [\G_{3'}^{\G_4},\G_5]=1 \\
Relation[45]: \ \ [\G_{3'}^{\G_{4'}},\G_{5'}]=1 \\
Relation[46]: \ \ [\G_2,\G_5]=1 \\
Relation[47]: \ \ [\G_2,\G_{5'}]=1 \\
Relation[48]: \ \ [\G_{2'},\G_5]=1 \\
Relation[49]: \ \ [\G_{2'},\G_{5'}]=1 \\
Relation[50]: \ \ <\G_1,\G_5>=1 \\
Relation[51]: \ \ <\G_1,\G_{5'}>=1 \\
Relation[52]: \ \ <\G_{1'},\G_5>=1 \\
Relation[53]: \ \ <\G_{1'},\G_{5'}>=1 \\ 
Relation[54]: \ \ [\G_3^{\G_4},
\G_{5'}^{(\G_6 \cdot \G_{6'} \cdot \G_5)}]=1 \\
Relation[55]: \ \ [\G_{3'}^{\G_4},
\G_{5'}^{(\G_6 \cdot \G_{6'} \cdot \G_5)}]=1 \\
Relation[56]: \ \ [\G_3^{(\G_4 \cdot \G_{4'} \cdot \G_4)},
\G_5^{(\G_6 \cdot \G_{6'} \cdot \G_5 \cdot \G_{5'} \cdot \G_5)}]=1 \\
Relation[57]: \ \ [\G_{3'}^{(\G_4 \cdot \G_{4'} \cdot \G_4)},
\G_5^{(\G_6 \cdot \G_{6'} \cdot \G_5 \cdot \G_{5'} \cdot \G_5)}]=1 \\
Relation[58]: \ \ [\G_2^{(\G_{1'} \cdot \G_1)},
\G_5]=1 \\
Relation[59]: \ \ [\G_2^{(\G_{1'} \cdot \G_1)},\G_{5'}]=1 \\
Relation[60]: \ \ [\G_{2'}^{(\G_{1'} \cdot \G_1)},
\G_5]=1 \\
Relation[61]: \ \ [\G_{2'}^{(\G_{1'} \cdot \G_1)},\G_{5'}]=1 \\
Relation[62]: \ \ \G_5^{(\G_1 \cdot \G_{1'} \cdot \G_2 \cdot 
\G_{2'} \cdot \G_5)}=\G_{5'}^{(\G_6 \cdot \G_{6'})} \\
Relation[63]: \ \ \G_1 \cdot \G_1=1 \\
Relation[64]: \ \ \G_{1'} \cdot \G_{1'}=1 \\
Relation[65]: \ \ \G_2 \cdot \G_2=1 \\
Relation[66]: \ \ \G_{2'} \cdot \G_{2'}=1 \\
Relation[67]: \ \ \G_3 \cdot \G_3=1 \\
Relation[68]: \ \ \G_{3'} \cdot \G_{3'}=1 \\
Relation[69]: \ \ \G_4 \cdot \G_4=1 \\
Relation[70]: \ \ \G_{4'} \cdot \G_{4'}=1 \\
Relation[71]: \ \ \G_5 \cdot \G_5=1 \\
Relation[72]: \ \ \G_{5'} \cdot \G_{5'}=1 \\
Relation[73]: \ \ \G_6 \cdot \G_6=1 \\
Relation[74]: \ \ \G_{6'} \cdot \G_{6'}=1.
$

\subsection{Relations $\hat{F}_1(\hat{F}_1)^{\rho^{-1}}$}
\indent

$(\hat{F}_1(\hat{F}_1)^{\rho^{-1}})^{\underline{Z}^{-2}_{5,66'}Z^2_{33',4}}$
are conjugated by $\underline{Z}^{-2}_{5,66'}Z^2_{33',4}$. Therefore, 
the corresponding elements $A, B$ in this case are conjugations of 
the corresponding $\G_j, \G_{j'}$ by 
$\underline{Z}^{-2}_{5,66'}Z^2_{33',4}, 1 \leq j \leq 6$.\\
We denote the conjugated  $\G_j$ and $\G_{j'}$ by 
$\bar\G_j$ and $\bar\G_{j'}$ respectively.

$\bar \G_j = \G_j$
for $j = 1, 1', 2, 2', 4'$.
 
We can view the other ones in the following figures.

\underline{$\bar \G_3,\bar \G_{3'},\bar \G_4 $}: \  $\vcenter{\hbox{\epsfbox{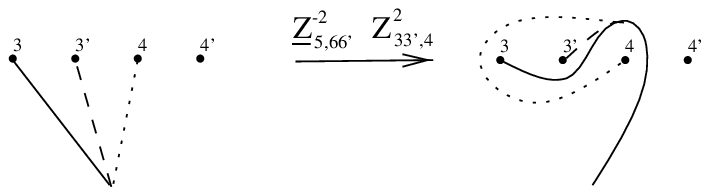}}}$ \\
$\bar\G_3 = \G_3^{\G_4} \\
\bar\G_{3'} = \G_{3'}^{\G_4} \\
\bar\G_4 = \G_4^{(\G_3 \cdot \G_{3'} \cdot \G_4)}$.\\
\\
\underline{$\bar \G_5$}: \  $\vcenter{\hbox{\epsfbox{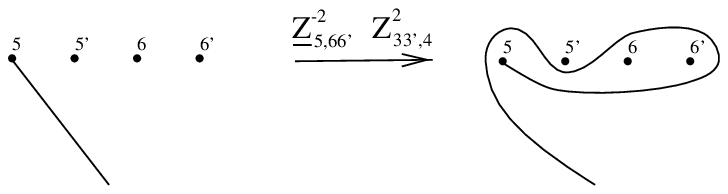}}}$ \\
$\bar\G_5 = \G_5^{(\G_{6'} \cdot \G_6 \cdot \G_5)}$.\\
\\
\underline{$\bar \G_{5'}$}: \  $\vcenter{\hbox{\epsfbox{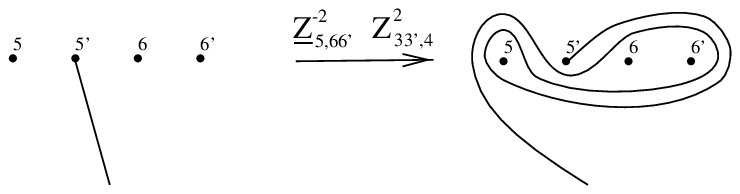}}}$ \\
$\bar\G_{5'} = \G_{5'}^{(\G_6 \cdot \G_{6'} \cdot \G_5 \cdot \G_{6'} \cdot \G_6 \cdot \G_5)}$.\\
\\
\underline{$\bar \G_6,\bar \G_{6'}$}: \  $\vcenter{\hbox{\epsfbox{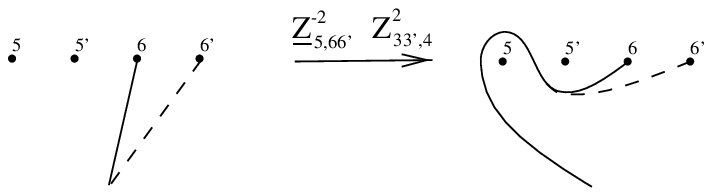}}}$ \\ 
$\bar\G_6 = \G_6^{\G_5} \\
\bar\G_{6'} = \G_{6'}^{\G_5}$.\\

The following relations are induced from the factors in
$(\hat{F}_1(\hat{F}_1)^{\rho^{-1}})^{\underline{Z}^{-2}_{5,66'}Z^2_{33',4}}$ 
and their paths in the second table of $H_{V_2}$.\\
$Relation[1]: \ \ <\bar\G_2,\bar\G_3>=1  \\
Relation[2]: \ \ <\bar\G_2,\bar\G_{3'}>=1  \\
Relation[3]: \ \ <\bar\G_{2'},\bar\G_3>=1  \\
Relation[4]: \ \ <\bar\G_{2'},\bar\G_{3'}>=1  \\
Relation[5]: \ \ <\bar\G_2,\bar\G_{3'}^{\bar\G_3}>=1  \\
Relation[6]: \ \ <\bar\G_{2'},\bar\G_{3'}^{\bar\G_3}>=1  \\
Relation[7]: \ \ [\bar\G_{3'},\bar\G_6]=1  \\
Relation[8]: \ \ [\bar\G_3^{\bar\G_{3'}},\bar\G_{6'}]=1  \\
Relation[9]: \ \ [\bar\G_3,\bar\G_{6'}^{\bar\G_6}]=1  \\
Relation[10]: \ \ [\bar\G_3^{(\bar\G_2 \cdot \bar\G_{2'} \cdot \bar\G_3)},
\bar\G_6]=1  \\
Relation[11]: \ \ [\bar\G_{3'}^{(\bar\G_2 \cdot \bar\G_{2'} \cdot \bar\G_{3'})},
\bar\G_{6'}]=1  \\
Relation[12]: \ \ [\bar\G_{3'}^{(\bar\G_3 \cdot \bar\G_2 \cdot \bar\G_{2'} \cdot \bar\G_3 \cdot \bar\G_{3'} \cdot \bar\G_3)},\bar\G_{6'}^{\bar\G_6}]=1  \\
Relation[13]: \ \ <\bar\G_2^{\bar\G_3},\bar\G_6>=1  \\
Relation[14]: \ \ <\bar\G_2^{\bar\G_{3'}},\bar\G_{6'}>=1  \\
Relation[15]: \ \ <\bar\G_{2'}^{\bar\G_3},\bar\G_6>=1  \\
Relation[16]: \ \ <\bar\G_{2'}^{\bar\G_{3'}},\bar\G_{6'}>=1  \\
Relation[17]: \ \ <\bar\G_2^{(\bar\G_3 \cdot \bar\G_{3'} \cdot \bar\G_3)},
\bar\G_{6'}^{\bar\G_6}>=1  \\
Relation[18]: \ \ <\bar\G_{2'}^{(\bar\G_3 \cdot \bar\G_{3'} \cdot \bar\G_3)},
\bar\G_{6'}^{\bar\G_6}>=1 \\
Relation[19]: \ \ \bar\G_1=\bar\G_{2'}^{(\bar\G_3 \cdot \bar\G_6)} \\
Relation[20]: \ \ \bar\G_{1'}=\bar\G_2^{(\bar\G_{2'} \cdot \bar\G_3 \cdot 
\bar\G_6)} \\
Relation[21]: \ \ \bar\G_1=\bar\G_{2'}^{(\bar\G_3 \cdot \bar\G_{3'} \cdot \bar\G_3 
\cdot \bar\G_6 \cdot \bar\G_{6'} \cdot \bar\G_6)} \\
Relation[22]: \ \ \bar\G_{1'}=\bar\G_2^{(\bar\G_{2'} \cdot \bar\G_3 \cdot \bar\G_{3'} 
\cdot \bar\G_3 \cdot \bar\G_6 \cdot \bar\G_{6'} \cdot \bar\G_6)} \\  
$

\subsection{Relations $C_i$}
\indent

In the same way as above, we derive relations from the braids  $C_i$. The list appears in [4 , Section 4.11] in a global numeration of the generators.

\section{Results}
In this paper we presented the list of braids and the braid monodromy factorization $\Delta^2_{54}$ corresponding to the branch curve $S$ of $T \times T$.
Recall that $S$ compounds nine curves, each one of them is treated in [4 ].
The computations in detail appear in [4 , Chapter3].

We quoted the van Kampen Method and Theorems and we used them to obtain relations for the fundamental group $\pi_1(\C^2-S,M)$. We presented some relations, the complete list appears in [4 , Chapter4]. Setting $\G_j^2=\G_{j'}^2=1$, we can get a presentation for the group $\tilde\pi_1$. 

Through this paper we concentrated on the 6-point $V_2$ and showed complete computations concerning this point.As quoted above, all other computations appear in [4 , Chapters 3 and 4].

In [7] we compute the fundamental group $\pi_1((T \times T)_{Gal})$ 
of the Galois cover of $T \times T$ with respect to a generic projection to $\C\P^2$. Recall that the fundamental group  $\pi_1((T \times T)_{Gal}^{Aff})$ is the kernel of the surjection $\psi: \tilde\pi_1 \rightarrow S_{18}$.

The Galois cover is a surface of a general type. 
In [7] we verify the Bogomolov Conjecture.  Bogomolov conjectured that 
if a surface of a general type has a positive index, then it has an
infinite fundamental group. 
  
Let $X$ be a surface, $f : X \rightarrow \C\P^2$ is 
a generic projection, $S \subset \C\P^2$ its branch curve, 
$m = deg S, d = \# $ nodes in $S, \rho = \#$ of cusps in 
$S, \mu = \#$ of tangency points in $S$ for 
a generic projection to $\C\P^1, n = deg f$. 
Then $C^2_1 (X_{Gal}) = \frac{ n!}{4} (m-6)^2$,  $C_2 (X_{Gal}) = 
n! (\frac{ m^2}{2} - \frac{ 3m}{2} + 3 - \frac{ 3d}{4} - \frac{ 4 \rho}{3})$ 
and the index $\tau=\frac{1}{3}(C^2_1(X_{Gal})-2C_2(X_{Gal}))$,  see 
[16, pp. 603-604].

We compute the index of  $(T \times T)_{Gal}$.\\
By the above computations, 
$n = 18, \ \mu = 54, \ \rho = 216,\  d = 1080, \ m = 54$.  
$C^2_1  (X_{Gal}) = 576 \cdot 18!$ and $C_2  (X_{Gal}) = 282 \cdot
18!$.  Therefore $\tau(X_{Gal}) = 4 \cdot 18!$.
The index is positive and the group we derive in [7] is infinite.

\section{Notations}\label{not}
\begin{list}{}{}
\item $(A)_B = B^{-1} AB = A^B.$
\item $X$ an algebraic surface, $X \subseteq \C \P^n$.
\item $X_0$ a degenerated object of a surface $X \ , \ X_0 \subset \C \P^N$.
\item $S$ an algebraic curve defined over $\R, \ S \subset \C^2$.
\item$E(\mbox{ resp}.D)$ be a closed disk on the $x$-axis (resp. y-axis) with the center on the real part of the $x$-axis (resp. y-axis), \st
$\{$singularities of $\pi_1\} \subseteq E \ast (D-\partial D)$.
\item $\pi  : S \rightarrow E$.
\item  $K(x) = \pi^{-1}(x)$.
\item $N = \{x \in E \suchthat \#K(x) < n \}$.
\item $u$ real number \st  $x << u \;\;\;  \forall  x \in N$.
\item $\C_u = \pi^{-1}(u)$.
\item $\rho_j = Z_{jj'}$.
\item $\vp$ = the braid monodromy of an algebraic curve $S$ in $M$.
\item $B_p[D,K]$ = the braid group.
\item $\Delta^2_p = (H_1 \cdots H_{p-1})^p$.
\item $\pi_1(\C^2 - S, M)$ = the fundamental group of a complement of a
branch curve $S$.
\item $\tP_1 = \frac{\pi_1(\C^2 - S, M)}{\la \G^2_j, \G^2_{j'}\ra}$ \ .
\item $T$ = complex torus.
\item $\underline{z}_{ij}$(resp.  $\bar{z}_{ij}$) = a path from $q_i 
\mbox{ to } q_j$ below(resp. above) the real line.
\end{list}
The corresponding halftwists are:
$
H(\underline{ z}_{ij}) = \underline{Z}_{ij} \; ; \; H(\bar{z}_{ij}) = \bar{Z}_{ij}.$

\end{document}